\documentclass[12pt,a4paper]{article}
\pdfoutput=1

\usepackage{amsmath}
\usepackage{amssymb}
\usepackage{latexsym}
\usepackage{bbm}
\usepackage{graphicx}

\begin{document}

\newtheorem{pb}{Problem}
\newtheorem{thm}{Theorem}
\newtheorem{prop}{Proposition}
\newtheorem{cor}{Corollary}
\newtheorem{lem}{Lemma}
\newtheorem{rmq}{Remark}
\newtheorem{dfn}{Definition}

\newcommand{\mc}[1]{\mathcal{#1}}
\newcommand{\mb}[1]{\mathbf{#1}}
\newcommand{\mbb}[1]{\mathbb{#1}}
\newcommand{\mbbm}[1]{\mathbbm{#1}}
\newcommand{\be}{\begin{equation}}
\newcommand{\ee}{\end{equation}}
\newcommand{\bea}{\begin{eqnarray}}
\newcommand{\eea}{\end{eqnarray}}
\newcommand{\ba}{\begin{array}}
\newcommand{\ea}{\end{array}}
\newcommand{\baa}{\left\{ \begin{array}}
\newcommand{\eaa}{\end{array} \right.}
\newcommand{\dt}{\partial_t}
\newcommand{\dx}{\partial_x}
\newcommand{\dy}{\partial_y}
\newcommand{\dz}{\partial_z}
\newcommand{\dzz}{\partial_{zz}}
\newcommand{\dxx}{\partial_{xx}}
\newcommand{\dyy}{\partial_{yy}}
\newcommand{\dxya}{\partial_{xy}^{\alpha}}
\newcommand{\dint}{\displaystyle \int}
\newcommand{\diint}{\displaystyle \iint}
\newcommand{\Frac}{\displaystyle \frac}
\newcommand{\dsum}{\displaystyle \sum}
\newcommand{\tld}[1]{\widetilde{#1}}
\newcommand{\udl}[1]{\underline{#1}}
\newcommand{\ovl}[1]{\overline{#1}}
\newcommand{\bproof}{\paragraph{Proof. }}
\newcommand{\eproof}{$\Box$}
\newcommand{\Ker}{\textrm{Ker\,}}
\newcommand{\Ima}{\textrm{Im\,}}
\newcommand{\Arg}{\textrm{Arg\,}}
\newcommand{\hookuparrow}{\cup\hskip-4.15pt{ }^{{}^{\textstyle\uparrow}}}
\newcommand{\etal}{{\it et al.}}
\newcommand{\divg}{\textrm{div}}
\newcommand{\Supp}{\textrm{Support }}
\newcommand{\loc}{\textrm{loc}}
\newcommand{\la}{\langle}
\newcommand{\ra}{\rangle}


\title{Regularity results for the Primitive Equations of the ocean}

\author{Ma\"elle Nodet}


\maketitle

\noindent Affiliation: Laboratoire de Mod\'elisation et de Calcul, Universit\'e Joseph Fourier Grenoble 1, CNRS\\
Address: LMC-IMAG, 51 rue des Math\'ematiques, BP53, 38041 Grenoble cedex 9, France.\\
E-mail: maelle.nodet@imag.fr\\

\paragraph{Abstract.} We consider the linear Primitive Equations of the ocean in the three dimensional space, with horizontal periodic and vertical Dirichlet boundary conditions. Thanks to Fourier transforms we are able to calculate explicitly the pressure term. We then state existence, unicity and regularity results for the linear time-depending Primitive Equations, with low-regularity right-hand side. 

\section{Introduction and main results}

We establish regularity results for the linear Primitive Equations (PE) of the ocean in the three dimensional space. For the nonlinear PE, the first work of Lions, Temam and Wang \cite{LionsTemam92} and the further paper by Temam and Ziane \cite{TemamZiane04} proved global existence of weak solutions and local existence and unicity of strong solutions. The regularity of the linear Stokes-type problem related to the PE has been studied by Ziane \cite{Ziane95}, Hu, Temam and Ziane \cite{TemamHu02} and Temam and Ziane \cite{TemamZiane04}. For diverse boundary conditions, the authors prove the regularity of weak solutions, when the right-hand side stays in $L^2$.\\
Our work is motivated by the following remark: many problems involving the PE of the ocean (such as numerical ocean modelling, assimilation of surface data or more theoretically controllability of the PE) need the calculation, or at least the estimation, of the pressure term. But in the previous studies of the PE, the regularity of the pressure is not explicitly investigated. Thus the aim of this paper is to calculate explicitly the most singular part of the pressure term in order to obtain more precise regularity results, in particular with less regular right-hand side. The analogous of this question for the full Stokes problem has been addressed by Fabre and Lebeau \cite{Fabre02}. \\
The paper is organized as follows: in the rest of this section our main results are stated, in section 2 we present some preliminaries, in section 3 we prove theorem \ref{thm:1}, some further remarks and results are provided in section 4.

\subsection{The Primitive Equations of the ocean}
The Primitive Equations of the ocean are at the base of general ocean circulation models, intensively used by oceanographers. The linear equations we are interested with are the following:
\be 
\baa{ll} 
\label{eq:1}
\dt u - \nu \Delta u  - \alpha\, v +  \dx p   =  f_1 &\textrm{ in } \Omega \times (0,T)\\
\dt v -  \nu \Delta v  + \alpha\, u +  \dy p  =  f_2 &\\
 \dz p  -  \beta \,\theta =0  &\\ 
\partial_{t} \theta - \nu \Delta \theta   +\gamma\, w =  f_3  &\medbreak \\
w(x,y,z,t)=-\int_0^z \dx u(x,y,z',t) + \dy v(x,y,z',t) \, dz'& 
\eaa  
\ee
with the following initial conditions:
\be
\label{eq:3}
U(t=0) = U_0, \qquad \theta(t=0) = \theta_0 \textrm{ in } \Omega
\ee
and boundary conditions:
\be 
\baa{l} 
\label{eq:2} 
u,v,w,\theta,p  \textrm{ are periodic in } x,y \medbreak\\\medbreak
u=0,v=0,\theta=0 \quad \textrm{ on } \mbb{T}^2\times \{z=0,z=a\}\times(0,T) \\
\int_{z=0}^a \dx u + \dy v \, dz= 0 \quad \textrm{ on } \mbb{T}^2\times(0,T)
\eaa 
\ee
where\\
\indent-- $\Omega$ is a horizontally periodic and vertically bounded ocean basin: $\Omega = \mbb{T}^2 \times(0,a)$, with $\mbb{T}^2=(\mbb{R}/2\pi\mbb{Z})^2$ the bidimensional torus;\\
\indent--  $U=(u,v)$ is the horizontal velocity of the fluid, $w$ its vertical velocity;\\
 \indent--  $\theta$ is the temperature around a vertical temperature profile $\theta=\tld{\theta}-\theta_b - z\frac{\theta_b-\theta_a}{a}$, with $\tld{\theta}$ the real temperature, $\theta_a$ and $\theta_b$ top and bottom boundary conditions for $\tld{\theta}$; \\
 \indent-- $p$ is the pressure;\\
 \indent--   $F=(f_1,f_2,f_3)$ is a given forcing term;\\
 \indent--  $\alpha$ is the constant Coriolis parameter;\\
 \indent--  $\nu$ is the kinematic viscosity and the temperature diffusion parameter;\\
 \indent--  $\beta$ is a physical constant, depending on the gravity constant;\\
 \indent--  $\gamma$ is a constant.
\begin{rmq}
\begin{enumerate}
\item We have assumed, without loss of generality, that the salinity does not appear in the state equation, so that it is a passive tracer (see \cite{TemamZiane04} for the full equations).
\item In order to lighten notations, we have assumed that the kinematic horizontal and vertical viscosity and the diffusion parameter in the temperature equations are equal, we thus have similar results with different values of these parameters.
\item We choose to use Dirichlet boundary conditions because these are realistic physical conditions; moreover this enables us to consider low regularity forcings $F$.
\end{enumerate}
\end{rmq}
We will also use the stationary linear model $\mc{S}_\lambda$ for spectral study, with  $\lambda \in \mbb{C}$:
\be 
\ba{c}
\mc{S}_\lambda(u,v,\theta) = F \smallbreak\\
\Updownarrow  \smallbreak \\
\baa{ll} 
\label{eq:4}
\lambda u - \nu \Delta u  - \alpha\, v +  \dx p   =  f_1 &\textrm{ in } \Omega \\
\lambda v -  \nu \Delta v  + \alpha\, u +  \dy p  =  f_2 &\\
 \dz p  -  \beta \,\theta =0  &\\ 
\lambda \theta - \nu \Delta \theta   +\gamma\, w =  f_3  & \medbreak\\
w(x,y,z)=-\int_0^z \dx u(x,y,z') + \dy v(x,y,z') \, dz'& 
\eaa
\ea
\ee
with stationary boundary conditions: 
\be 
\baa{l} 
\label{eq:5} 
u,v,w,\theta,p  \textrm{ are periodic in } x,y \medbreak\\\medbreak
u=0,v=0,\theta=0 \quad \textrm{ on } \mbb{T}^2\times \{z=0,z=a\} \\
\int_{z=0}^a \dx u + \dy v \, dz= 0 \quad \textrm{ on } \mbb{T}^2
\eaa 
\ee
\subsection{Some functional spaces}
Let us now introduce some functional spaces
\begin{dfn}\label{dfn:1}
Forall  $s\in\mbb{R}$
\be
\label{eq:6}
\ba{rl}
\mc{H}^{s}(\Omega) = \big\{ & f(x,y,z) = \sum_{k\in\mbb{N}^*,\zeta\in\mbb{Z}^2} f_{k,\zeta} \,e_k(z) \,e_\zeta(x,y),  (x,y,z)\in\Omega,\smallbreak\\
&\sum_{k\in\mbb{N}^*,\zeta\in\mbb{Z}^2} (1+\nu k^2+\nu|\zeta|^2)^{s} |f_{k,\zeta}|^2 < \infty \big\}
\ea
\ee
with:
\be
\label{eq:7}
\ba{rcll}
e_k(z)&=&\sqrt{\frac{2}{a}} \sin(\frac{k\pi z}{a})&\forall z\in(0,a)\smallbreak \\
e_\zeta(x,y)&=&\frac{1}{2\pi}e^{i(\xi x +\eta y)} &\forall (x,y)\in\mbb{T}^2
\ea
\ee
forall $k\in\mbb{N}^*$ and $\zeta=(\xi,\eta)\in\mbb{Z}^2$.
\end{dfn}
For $f\in\mc{H}^s$ we denote
\be
\label{eq:8}
\|f\|^2_{s} = \sum_{k\in\mbb{N}^*,\zeta\in\mbb{Z}^2} (1+\nu k^2+\nu|\zeta|^2)^{s} |f_{k,\zeta}|^2
\ee
$\mc{H}^{s}(\Omega)$ is a Hilbert space with the following inner product:
\be
\label{eq:9}
\langle f , g \rangle_{s} = \sum_{k\in\mbb{N}^*,\zeta\in\mbb{Z}^2} (1+\nu k^2+\nu|\zeta|^2)^{s} f_{k,\zeta}\, \ovl{g_{k,\zeta}}.
\ee
The following characterization holds true, where $H^s(\Omega)$ denotes the usual Sobolev space:
\begin{lem}\label{lem:1} 
\be
\label{eq:10}
\ba{rclcl}
-\frac32 < s< \frac12&\Rightarrow &\mc{H}^{s}(\Omega) &= &H^s(\Omega)\smallbreak \\
\frac12 < s <\frac52&\Rightarrow  &\mc{H}^{s}(\Omega) &= &\{ f\in H^s(\Omega), f|_{z=0}=f|_{z=a}=0 \}
\ea
\ee
\end{lem}
We define also the following classical spaces (see \cite{Girault79}, \cite{Temam84} or \cite{TemamHu02}):
\begin{dfn}\label{dfn:2} Let
\be
\label{7eq:21}
\ba{rcl}
E_1 &=& \{ U=(u,v)\in\mc{C}^\infty(\Omega)^2, u,v\textrm{ periodic in }x,y, \\
&&\quad u=0,v=0 \textrm{ on } \mbb{T}^2\times \{z=0,z=a\}\\
&&\quad\int_{0}^a \dx u(x,y,z') +\dy v(x,y,z') \, dz' =0, \forall (x,y)\in \mbb{T}^2 \}\medbreak\\
E_2 &=& \{ \theta\in\mc{C}^\infty(\Omega),  \theta \textrm{ periodic in } x,y, \\
&&\quad\theta=0 \textrm{ on } \mbb{T}^2\times \{z=0,z=a\}\}
\ea
\ee
Then $\mc{H}_1$ (respectively $\mc{H}_2$) is defined to be the closure of $E_1$ in $L^2(\Omega)^2$ (resp. $L^2(\Omega)$),  and $\mc{V}_1$ (resp. $\mc{V}_2$) is the closure of $E_1$ (resp. $E_2$) in $H^1(\Omega)^2$ (resp. $H^1(\Omega)$), and finally $\mc{H}=\mc{H}_1\times \mc{H}_2$,  $\mc{V}=\mc{V}_1\times \mc{V}_2$.
\end{dfn}
Inner products on $\mc{H}$ and $\mc{V}$ are:
\be
\ba{rcl}
\label{5eq:364}
(X,X')_{\mc{H}} &=& (u,u')_{L^2(\Omega)} + (v,v')_{L^2(\Omega)} + \frac{\beta}{\gamma} (\theta,\theta')_{L^2(\Omega)}\\
&=& \int_{\Omega} ( u \ovl{u'} + v \ovl{v'}+\frac{\beta}{\gamma} \theta \ovl{\theta'}) \, dx\,dy\,dz \medbreak \\
(X,X')_{\mc{V}} &=& (u,u')_{H^1_0(\Omega)} + (v,v')_{H^1_0(\Omega)} + \frac{\beta}{\gamma} (\theta,\theta')_{H^1_0(\Omega)}\\
&=& \int_{\Omega} ( \nabla u .\nabla \ovl{u'} + \nabla v .\nabla \ovl{v'}+ \frac{\beta}{\gamma} \nabla \theta .\nabla \ovl{\theta'}) \, dx\,dy\,dz
\ea
\ee
\subsection{Results}
Our main result states as follows:
\begin{thm}\label{thm:1}
Let $\sigma \in ]-\frac32,\frac12[$, $\sigma\neq -\frac12$.\\
Let $F(t) = (f_1,f_2,f_3) \in (L^2(\mbb{R} ;\mc{H}^\sigma)^3)$ with $\Supp(F) \subset \{t\geq 0\}$.\\
There exists a unique
\be
\label{eq:11}
X(t)=(u,v,\theta)\in (L^2(\mbb{R} ; \mc{H}^{\sigma+2})^3), \quad \Supp(X) \subset \{t\geq 0\}
\ee
and there exists a unique (up to a distribution depending only on $t$) pressure
\be
\label{eq:12}
p(t)\in\mc{D}'(\mbb{R}\times\Omega), \quad \Supp(p) \subset \{t\geq 0\}
\ee
so that the following equation holds true, in the sense of distributions in $\mbb{R}\times\Omega$
\be
\label{eq:13}
\ba{l}
\baa{l}
\dt u -\nu \Delta u -\alpha v +\dx p = f_1\\
\dt v -\nu \Delta v + \alpha u + \dy p = f_2\\
\dz p - \beta \theta = 0\\
\dt \theta -\nu\Delta\theta+\gamma w =f_3 
\eaa \smallbreak\\
\textrm{ with }w(z)=-\int_0^z (\dx u +\dy v) \textrm{ and } w(a)=0
\ea
\ee
Moreover
\be
\label{eq:14}
\|X\|_{(L^2(\mbb{R} ; \mc{H}^{\sigma+2}))^3} \,\leq\, C\, \|F\|_{(L^2(\mbb{R} ;\mc{H}^\sigma))^3}
\ee
and the temperature $\theta$ verifies
\be
\label{eq:15}
\dt \theta \in L^2(\mbb{R} ;\mc{H}^\sigma)\quad\textrm{ and }\quad \|\dt \theta\|_{L^2(\mbb{R} ; \mc{H}^{\sigma})} \,\leq\, C\, \|F\|_{(L^2(\mbb{R} ;\mc{H}^\sigma))^3}
\ee
The pressure $p$ verifies
\be
\label{eq:16}
p(t,x,y,z) = c(t) + q(t,x,y) + \beta\int_0^z \theta(t,x,y,z')\,dz'
\ee
with 
\be
\label{eq:17}
c(t) \in\mc{D}'(\mbb{R}), \quad\Supp(c) \subset \{ t\geq 0 \}
\ee
and\\
\indent - for $\sigma \in ]-\frac12, \frac12[$ we have
\be
\label{eq:18}
q(t,x,y) \in L^2(\mbb{R} ; H^{\sigma+1}(\mbb{T}^2))\quad\textrm{ and }\quad\|q \|_{L^2(\mbb{R} ; H^{\sigma+1}(\mbb{T}^2))} \,\leq\,C\,\|F\|_{(L^2(\mbb{R} ;\mc{H}^\sigma))^3}\ee
\indent - for $\sigma \in ]-\frac32, -\frac12[$ we have
\be
\label{eq:19}
\ba{c}
q(t,x,y) = q_1(t,x,y)+ q_2(t,x,y)\medbreak\\
q_2(t,x,y) \in L^2(\mbb{R} ; H^{\sigma+1}(\mbb{T}^2))\textrm{ and }\|q_2 \|_{L^2(\mbb{R} ; H^{\sigma+1}(\mbb{T}^2))} \,\leq\,C\,\|F\|_{(L^2(\mbb{R} ;\mc{H}^\sigma))^3} \smallbreak\\
q_1(t,x,y) \in H^{\sigma/2+1/4}(\mbb{R} ; H^{1}(\mbb{T}^2))\textrm{ and }\|q_1\|_{H^{\sigma/2+1/4}(\mbb{R} ; H^{1}(\mbb{T}^2))} \,\leq\,C\,\|F\|_{(L^2(\mbb{R} ;\mc{H}^\sigma))^3}
\ea
\ee
\end{thm}
\begin{rmq}\label{rmq:2}
\begin{enumerate}
\item The regularity exponent $\sigma$.  ~\\[0.1cm]
- With a forcing term $F\in(L^2(\mbb{R} ;\mc{H}^\sigma)^3)$, we cannot have more than  $X\in (L^2(\mbb{R} ; \mc{H}^{\sigma+2})^3)$. Thus the boundary condition  $X|_{z=0,z=a}=0$ is well defined only if  $\sigma+2>\frac12$, i.e. $\sigma >-\frac32$.\\
-  $\sigma=-\frac12$ is a critical exponent for the regularity of the pressure, whose description is more technical.\\
- Considering only $\sigma <\frac12$ enables us to use the spaces  $\mc{H}^\sigma$ and to do explicit calculations.
\item An explicit formula for $q_1$. ~\\[0.1cm]
Actually, for $\sigma\in]-\frac32,-\frac12[$, we will prove a more precise result than (\ref{eq:19}), namely
\be
\label{eq:20}
q(t,x,y) = q_1(t,x,y)+ q_2(t,x,y), \quad q_2(t,x,y) \in L^2(\mbb{R} ; H^{\sigma+1}(\mbb{T}^2))
\ee
where $q_1$ is explicit as a function of $F$ (see remark \ref{5rmq:12} and formula (\ref{5eq:228})).
\item Formula (\ref{eq:16}). ~\\[0.1cm]
In formula (\ref{eq:16}) for $p$, $q$ is the value of $p$ at $z=0$. We can replace $q$ either by $p(t,x,y,z_0)$, for any $z_0\in[0,a]$, or by $\int_0^a p(t,x,y,z)\,dz$, the results remain the same.
\item Maximal estimates. ~\\[0.1cm]
For $\sigma >-\frac12$, we have $\dx p,\dy p \in L^2(\mbb{R} ;\mc{H}^\sigma)$, so that the pressure gradient term can be seen as a forcing term and we have the following maximal estimates
\be
\label{eq:21}
\|X\|_{(L^2(\mbb{R} ; \mc{H}^{\sigma+2}))^3} + \|\dt X\|_{(L^2(\mbb{R} ; \mc{H}^{\sigma}))^3}\,\leq\, C\, \|F\|_{(L^2(\mbb{R} ;\mc{H}^\sigma))^3}
\ee
However, for $\sigma\in]-\frac32,-\frac12[$, the maximal estimate is wrong (see remark \ref{5rmq:13}).
\end{enumerate}
\end{rmq}
We prove also the following corollary for the Cauchy problem with $\sigma=-1$:
\begin{cor}\label{cor:1}
Let $\varphi(t) \in \mc{C}_c^{\infty}(]0,T[)$, $F(t) = (f_1,f_2,f_3) \in (L^2(0,T;\mc{H}^{-1}))^3$, $X_0 \in \mc{H}$. Let $(X,p)$ be the unique solution of equation (\ref{eq:13}) with
\be
\label{eq:22}
\ba{l}
X=(u,v,\theta) \in L^2(0,T;\mc{V}) \cap \mc{C}([0,T];\mc{H}),\quad X(t=0)=X_0 \\
p\in \mc{D'}(0,T;L^2(\Omega))
\ea
\ee
Then $\varphi p$ is rewritten as
\be
\label{eq:23}
\varphi p(t,x,y,z) = c(t) + q(t,x,y) + \beta\int_0^z \theta(t,x,y,z')\,dz',\textrm{ with }c(t) \in\mc{D}'(\mbb{R})
\ee
with $q(t) \in H^{-1/4}(0,T;H^1(\mbb{T}^2))$ and we have
\be
\label{eq:24}
\ba{c}
q(t) \in L^2(0,T;L^2(\mbb{T}^2))\smallbreak\\
\Updownarrow\smallbreak\\
  \Delta_2^{-1} \big[ \int_0^a (\dt-\nu\Delta)^{-1}[\varphi\dx f_1+\varphi\dy f_2]\,dz \big] \in L^2(0,T;L^2(\mbb{T}^2))
\ea
\ee
where $\Delta_2$ is the horizontal Laplacian operator, defined by $\Delta_2 \psi = \dxx \psi + \dyy \psi$.
\end{cor}
\begin{rmq} \label{rmq:cor}
The preceding Cauchy problem can easily be addressed thanks to classical variational methods, but it gives less precise results regarding the pressure, see lemma \ref{lem:5} and remark \ref{5rmq:9}.
\end{rmq}
\section{Preliminary results}
%
%
%
%
%
\subsection{The Primitive Equations operator}
\label{5sec:2}
Multiplying equation (\ref{eq:4}) by $\gamma \ovl{u'}$, $\gamma \ovl{v'}$, $\gamma \ovl{w'}$, $\beta \ovl{\theta'}$ (with $X'=(u',v',\theta') \in \mc{V}$) and integrating by parts (using boundary conditions (\ref{eq:2})), we obtain formally:
\be
\label{5eq:52}
\ba{c}
\mc{S}_{\lambda}(X)=F \\ \Updownarrow  \\
\lambda (X,X')_{\mc{H}} + \nu (X,X')_{\mc{V}} + \beta B(X,X') + \alpha C(X,X') = (F,X')_{\mc{H}}, \quad \forall X' \in \mc{V}
\ea
\ee
where $B$ and $C$ are given by:
\be
\label{5eq:53}
\ba{lll}
B(X,X') = -  (\theta,w')_{L^2(\Omega)} + (w,\theta')_{L^2(\Omega)}\\
\quad \textrm{ with }w=-\int_0^z \dx u + \dy v, \quad  w'=-\int_0^z \dx u' + \dy v'\medbreak\\
C(X,X') = -  (v,u')_{L^2(\Omega)} +  (u,v')_{L^2(\Omega)}
\ea
\ee
and
\be
\label{5eq:54}
\ba{llll}
B(X,X) &= \int_{\Omega}( - \theta \ovl{w} + w \ovl{\theta}) &= 2i \Im( \int_{\Omega}w \ovl{\theta}) &\in i\mbb{R}\\
C(X,X)&= \int_{\Omega}( -v \ovl{u} + u \ovl{v}) &= 2 i \Im(\int_{\Omega} u \ovl{v} ) &\in i\mbb{R}
\ea
\ee
We define then $A(X,X') = (X,X')_{\mc{V}}$. The operator  $P=\nu A+\beta B+\alpha C$, called Primitive Equations operator,  maps $\mc{V}$ to $\mc{V'}$, and for all $(X,X') \in \mc{V}$ we have
\be
\ba{l}
\label{5eq:57}
\langle P(X),X'\rangle _{\mc{V'},\mc{V}} \, =  \,\nu (X,X')_{\mc{V}} + \beta B(X,X') + \alpha C(X,X')\medbreak\\
\langle (\lambda  + P)(X),X'\rangle _{\mc{V'},\mc{V}} \, =  \, \lambda (X,X')_{\mc{H}}+ \nu (X,X')_{\mc{V}} + \beta B(X,X') + \alpha C(X,X')
\ea
\ee
\begin{rmq} \label{5rmq:8}
Operator $A$ corresponds to the uncoupled Stokes-type equation obtained from  (\ref{eq:4}) with $\alpha=\beta=\gamma=0$, $B$ corresponds to the coupling (via the parameters $\beta$ and $\gamma$) between the vertical velocity $w$ and the temperature $\theta$ and  $C$ is the Coriolis operator.
\end{rmq}
We have:
\begin{lem} \label{lem:2}
The mapping
$$\Phi\,:\,(X,X')\mapsto \langle P(X),X'\rangle_{\mc{V'},\mc{V}} $$
 is continuous on $\mc{V}^2$. More precisely
\be
\label{5eq:113}
|\langle  P(X),X'\rangle_{\mc{V'},\mc{V}} | \leq ( \nu + 2 \frac{a^2}{\pi}\sqrt{\beta \gamma } + \frac{2 \alpha a^2 }{\pi^2}) \|X\|_{\mc{V}}\|X'\|_{\mc{V}} 
\ee
\end{lem}
\bproof Let $X=(u,v,\theta)$ and $X'=(u',v',\theta')$ be in $\mc{V}$. We have clearly:
\be
\label{5eq:58}
\ba{rcl}
\langle A(X),X'\rangle _{\mc{V'},\mc{V}} &=& (X,X')_{\mc{V}}\,\leq\,\|X\|_{\mc{V}}\|X'\|_{\mc{V}} \medbreak\\
|B(X,X')|& \leq&\|\theta\|_{L^2(\Omega)} \|w'\|_{L^2(\Omega)} + \|w\|_{L^2(\Omega)} \|\theta'\|_{L^2(\Omega)}\medbreak\\
|C(X,X')|&\leq& \|v\|_{L^2(\Omega)} \|u'\|_{L^2(\Omega)} + \|u\|_{L^2(\Omega)} \|v'\|_{L^2(\Omega)}
\ea
\ee
For all $\varphi \in H^1_0(\Omega)$ the following Poincar\'e inequality holds:
\be
\label{5eq:59}
\|\varphi\|^2_{L^2(\Omega)} \leq \frac{a^2}{\pi^2} \|\nabla \varphi\|^2_{L^2(\Omega)}
\ee
Thus we obtain
\be
\label{5eq:60}
\|u\|_{L^2(\Omega)} \leq  \|X\|_{\mc{H}}\leq \frac{a}{\pi} \|X\|_{\mc{V}}, \qquad \|\theta\|_{L^2(\Omega)} \leq \sqrt{\frac{\gamma}{\beta}}\|X\|_{\mc{H}}\leq \frac{a}{\pi} \sqrt{\frac{\gamma}{\beta}}\|X\|_{\mc{V}}
\ee
Using Cauchy-Schwarz inequality, we obtain
\be
\label{5eq:61}
\ba{rcl}
\|w\|^2_{L^2(\Omega)}&=&\|\int_0^z \dx u + \dy v\|^2_{L^2(\Omega)}\\
&\leq&  a^2 \|X\|^2_{\mc{V}}
\ea
\ee
A straightforward calculation gives the desired conclusion.\\
\eproof
\subsection{Qualitative spectral study}
\begin{dfn}\label{dfn:3}
We call eigenvalue of  $-P:\mc{V}\rightarrow\mc{V'}$ a complex number $\lambda$ so that $\lambda + P$ is not injective. We denote $\mbb{V}_P  $ the set of the eigenvalues of $-P$:
\be
\label{5eq:166}
\mbb{V}_P\, = \, \{\lambda \in \mbb{C},\, \exists X \in \mc{V}, \,X\neq 0,\, \langle\lambda X + P(X),X'\rangle_{\mc{V'},\mc{V}}=0\, \forall X'\in\mc{V} \}
\ee
\end{dfn}
\begin{lem} \label{lem:3} We have the following inclusion:
\be
\label{5eq:87}
\mbb{V}_P \subset \Big\{ \,\lambda \in \mbb{C}, \Re(\lambda)  \leq -\frac{\nu\pi^2}{a^2}  \textrm{ and } 
|\Im(\lambda)|  \leq 2\alpha + 2 a\sqrt{\beta \gamma} \sqrt{-\frac{\Re(\lambda)}{\nu}} \,\Big\}
\ee
\end{lem}
\bproof If $\lambda$ is an eigenvalue of $-P$, then there exists $X \in \mc{V}, X \neq 0$ such that
\be
\label{5eq:82}
\lambda (X,X')_{\mc{H}} + \nu (X,X')_{\mc{V}} + \beta B(X,X') + \alpha C(X,X')= 0, \forall X' \in \mc{V}
\ee
Using (\ref{5eq:54}) and (\ref{5eq:82}) with $X'=X$, we have
\be
\label{5eq:84}
\ba{l}
\Re(\lambda) \|X\|^2_{\mc{H}} + \nu \|X\|^2_{\mc{V}}  = 0\medbreak\\
\Im(\lambda) \|X\|^2_{\mc{H}} + 2 \beta \Im (\int_{\Omega} w \ovl{\theta}) + 2\alpha \Im (\int_{\Omega} u \ovl{v}) = 0
\ea
\ee
Thanks to (\ref{5eq:60}) $\|X\|^2_{\mc{H}} \leq \frac{a^2}{\pi^2} \|X\|^2_{\mc{V}}$, we obtain:
\be
\label{5eq:85}
\Re(\lambda) = -\frac{\nu \|X\|^2_{\mc{V}}}{\|X\|^2_{\mc{H}}} \leq -\frac{\nu\pi^2}{a^2}
\ee
With (\ref{5eq:60}) and (\ref{5eq:61}) we have:
\be
\label{5eq:86}
\ba{rcl}
|\Im(\lambda)|  &=& 2 \beta |\Im (\int_{\Omega} w \ovl{\theta})| / \|X\|^2_{\mc{H}} + 2\alpha |\Im (\int_{\Omega} u \ovl{v})|/ \|X\|^2_{\mc{H}}\\
&\leq&  2 \beta \|w\|_{L^2(\Omega)}\|\theta\|_{L^2(\Omega)}  / \|X\|^2_{\mc{H}} + 2 \alpha \|u\|_{L^2(\Omega)}\|v\|_{L^2(\Omega)}/ \|X\|^2_{\mc{H}}\\
&\leq&    2 a\sqrt{\beta\gamma} \|X\|_{\mc{V}}  / \|X\|_{\mc{H}} + 2\alpha \\
&=& 2 \alpha + 2 a \sqrt{\beta\gamma} \sqrt{-\frac{\Re(\lambda)}{\nu}}
\ea
\ee
From (\ref{5eq:85}) and (\ref{5eq:86}) we get that (\ref{5eq:87}) holds true.\\
\eproof

\subsection{First existence and unicity results}
Let us finish this section by stating two lemmas, whose proofs are based on very classical use of the variational method, as in \cite{LionsTemam92} or in \cite{TemamZiane04}.
\begin{lem} \label{lem:4}
If $\lambda \in \mbb{C}\setminus\mbb{V}_P$ and $Y=(y_1,y_2,y_3) \in (H^{-1}(\Omega))^3$ then there exists a unique $X=(u,v,\theta) \in \mc{V}$ and there exists a pressure $p(x,y,z) \in L^2(\Omega)$, unique up to a constant, so that
\be 
\baa{ll} 
\label{5eq:64}
\lambda u - \nu \Delta u   -\alpha v+  \dx p   =  y_1 &\textrm{ in } \Omega \\
\lambda  v -  \nu \Delta v   + \alpha u +  \dy p =  y_2&\\
 \dz p  - \beta \theta  = 0 &\\ 
\lambda  \theta - \nu \Delta \theta    + \gamma w =  y_3  & \medbreak\\
w(x,y,z)=-\int_0^z \dx u(x,y,z') + \dy v(x,y,z') \, dz'& 
\eaa  
\ee
\end{lem}

\begin{lem} \label{lem:5}
Let $T >0$, $X_0 \in \mc{H}$ and $F=(f_1,f_2,f_3) \in L^2(0,T,(H^{-1}(\Omega))^3)$. Then there exists a unique 
\be
\label{5eq:75} 
X=(u,v,\theta) \in L^2(0,T;\mc{V}) \cap \mc{C}([0,T];\mc{H}) 
\ee
and there exists a pressure 
\be
\label{5eq:142}
p\in \mc{D'}(0,T;L^2(\Omega))
\ee
unique (up to a time distribution),  such that the following equation holds true in the sense of distributions in  $\Omega \times (0,T)$:
\be 
\baa{ll} 
\label{5eq:76}
\dt u - \nu \Delta u  -\alpha v +  \dx p   =  f_1   \\
\dt v -  \nu \Delta v +\alpha u  +  \dy p =  f_2 &\\
 \dz p  - \beta \theta  = 0 &\\ 
\dt \theta - \nu \Delta \theta   + \gamma w =  f_3   & \medbreak\\
\textrm{ with } w(z)=-\int_0^z \dx u(z') + \dy v(z') \, dz'
\eaa
\ee
and
\be
\label{5eq:170} 
(u,v,\theta)|_{t=0} = X_0 
\ee
\end{lem}
\begin{rmq} \label{5rmq:9}
The derivative $\frac{dX}{dt}$ is in $H^{-1}(0,T,(H^{1}_0(\Omega))^3)$ and equation (\ref{5eq:76}) tells us in particular:
\be
\label{5eq:178}
\nabla p \in L^2(0,T,(H^{-1}(\Omega))^3) + H^{-1}(0,T,(H^{1}_0(\Omega))^3)
\ee
\end{rmq}

\section{Proof of theorem \ref{thm:1}}

The proof is organized as follows. 
 In section \ref{thm1:1} we take the Fourier-Laplace transform of the equation, first in the horizontal coordinates, then in the vertical one and finally in time. Spectral parameters are then introduced, $\lambda=i\tau$ and $\omega=\lambda+\nu \zeta^2$, where $\tau$ is the Laplace parameter and $\zeta$ the horizontal Fourier variable. 
In section \ref{thm1:2} we introduce the function $M_\sigma(\lambda,\zeta)$ and we prove preliminary estimates for this function. 
In section \ref{thm1:3} we study the uncoupled system, ie (\ref{5eq:62}) ie with $\alpha=\beta=\gamma=0$. This is the core of the proof. We will use the function $M_\sigma$ in order to establish in theorem \ref{5prop:3} optimal estimates for the uncoupled system depending on the values and asymptotics of the parameters and the vertical Fourier variable. 
In section \ref{thm1:4} we use the results of section \ref{thm1:3} to establish estimates for the coupled system. 
We conclude the proof in section \ref{thm1:5}.

\subsection{First reductions}\label{thm1:1}

Estimation  (\ref{eq:15}) for temperature is straightforward: if (\ref{eq:14}) holds, then $w\in L^2(\mbb{R} ;\mc{H}^\sigma)$ and
\be
\label{5eq:33}
\|w\|_{L^2(\mbb{R} ;\mc{H}^\sigma)} \,\leq\,C\, \|F\|_{L^2(\mbb{R} ;\mc{H}^\sigma)}
\ee
then the temperature satisfies
\be
\label{5eq:34}
\dt \theta -\nu\Delta\theta = f_3 -\gamma w
\ee
and (\ref{eq:15})follows easily.
Thus it is sufficient to prove existence, unicity, (\ref{eq:16}), (\ref{eq:14}), (\ref{eq:18}) and (\ref{eq:19}).

\paragraph{Fourier transform in space.} 
\label{5sec:red1}

For $f\in\mc{D}'(\mbb{R}\times\Omega)$, we write
\be
\label{5eq:24}
f(t,x,y,z) = \sum_{\zeta\in\mbb{Z}^2} f_\zeta (t,z) \, e^{i\zeta.(x,y)}
\ee
Equation (\ref{eq:13}) is equivalent to the following equations, with parameter $\zeta=(\xi,\eta)\in\mbb{Z}^2$:
\be
\label{5eq:25}
\ba{l}
\baa{l}
\dt u_\zeta -\nu \dzz u_\zeta +\nu\zeta^2 u_\zeta-\alpha v_\zeta +i\xi p_\zeta = f_{1,\zeta}\\
\dt v_\zeta -\nu \dzz v_\zeta +\nu\zeta^2 v_\zeta+ \alpha u_\zeta + i\eta p_\zeta = f_{2,\zeta}\\
\dz p_\zeta - \beta \theta_\zeta = 0\\
\dt \theta_\zeta -\nu\dzz\theta_\zeta+\nu\zeta^2\theta_\zeta+\gamma w_\zeta =f_{3,\zeta} 
\eaa \smallbreak\\
\textrm{ with }w_\zeta(t,z)=-\int_0^z (i\xi u_\zeta +i\eta v_\zeta),\,\, w_\zeta(a)=0,\,\, X_\zeta|_{z=0,z=a}=0
\ea
\ee
The equation above gives
\be
\label{5eq:32}
p_\zeta(t,z) = p_\zeta(t,0) + \beta \int_0^z \theta_\zeta(t,z')\, dz'
\ee
So we set
\be
\label{5eq:274}
\ba{lclcl}
\zeta=0&: \qquad& c_0(t) = p_0(t,0) & ;& q_0(t) = 0\\
\zeta\neq 0 &:\qquad & c_\zeta(t) =  0 & ;& q_\zeta(t) =p_\zeta(t,0)
\ea
\ee

\paragraph{The spaces $H^s_\zeta$.}
We define now the following space of functions of  $z\in(0,a)$. For $f(z)=\sum_{k\in\mbb{N}^*}f_{k}\,e_k(z)$, we set
\be
\label{5eq:26}
\|f\|^2_{s,\zeta}=\sum_{k\in\mbb{N}^*} (1+\nu k^2 +\nu \zeta^2)^s |f_{k}|^2
\ee
and we denote by $H^s_\zeta$ the Hilbert space associated with this latter norm. Similarly to lemma \ref{lem:1} for spaces $\mc{H}^s$ we have 
\begin{lem}\label{5lem:4} 
\be
\label{5eq:27}
\ba{rclcl}
-\frac32 < s< \frac12&\Rightarrow & H^s_\zeta&= &H^s(0,a)\smallbreak \\
\frac12 < s <\frac52&\Rightarrow  & H^s_\zeta&= &\{ f(z)\in H^s(0,a), f|_{z=0}=f|_{z=a}=0 \}
\ea
\ee
\end{lem}
And for $f(t,x,y,z)=\sum_{\zeta\in\mbb{Z}^2}f_{\zeta}(t,z)\,e^{i\zeta.(x,y)}$  we obtain:
\be
\label{5eq:28}
\|f(t)\|^2_{\mc{H}^s} = \sum_{\zeta\in\mbb{Z}^2} \|f_\zeta(t,.)\|^2_{s,\zeta}
\ee
So that (\ref{eq:14}),  (\ref{eq:18}) and  (\ref{eq:19})  are equivalent to the following estimates, with $C$ independent of    $\zeta$:
\be
\label{5eq:29}
\ba{lcl}
\|X_\zeta\|_{(L^2(\mbb{R} ; H_\zeta^{\sigma+2}))^3} &\leq& C\, \|F_\zeta\|_{(L^2(\mbb{R} ;H_\zeta^\sigma))^3}\ea
\ee
and also for $q_\zeta = q_{1,\zeta}+q_{2,\zeta}$:
\be
\label{5eq:219}
\ba{lcl}
-\frac12 < \sigma < \frac12&:&\|q_\zeta\|_{L^2(\mbb{R}_+;H^{\sigma+1}_\zeta)} \, \leq\,  C \|F_\zeta\|_{L^2(\mbb{R}_+;H^{\sigma}_\zeta)}\\
-\frac32 < \sigma < -\frac12&:&\|q_{1,\zeta}\|_{H^{\sigma/2+1/4}(\mbb{R}_+;H^{1}_\zeta)}  + \|q_{2,\zeta}\|_{L^2(\mbb{R}_+;H^{\sigma+1}_\zeta)} \, \leq\,  C \|F_\zeta\|_{L^2(\mbb{R}_+;H^{\sigma}_\zeta)}
\ea
\ee

\paragraph{Case $\zeta=0$.} ~\\
In that case, $p_0$ vanishes from the first two equations of  (\ref{5eq:25}) , $w_0=0$ and (\ref{5eq:25}) gives
\be
\label{5eq:30}
\ba{l}
\baa{l}
\dt u_0 -\nu \dzz u_0 -\alpha v_0  = f_{1,0}\\
\dt v_0 -\nu \dzz v_0 + \alpha u_0 = f_{2,0}\\
\dt \theta_0 -\nu\dzz\theta_0  =f_{3,0} \\
p_0(t,z) = c_0(t) + \beta \int_0^z \theta_0(t,z')\, dz'
\eaa \smallbreak\\
\textrm{ with } X_0|_{z=0,z=a}=0
\ea
\ee
So that classical results on the heat equation give:
\be
\label{5eq:31}
\ba{lcl}
\|X_0\|_{(L^2(\mbb{R} ; H_0^{\sigma+2}))^3} &\leq& C\, \|F_0\|_{(L^2(\mbb{R} ;H_0^\sigma))^3}\smallbreak\\
\|\dt X_0\|_{(L^2(\mbb{R} ; H_0^{\sigma}))^3} &\leq& C\, \|F_0\|_{(L^2(\mbb{R} ;H_0^\sigma))^3}
\ea
\ee
Moreover $q_0(t)=0$ and estimates (\ref{eq:18}) and (\ref{eq:19}) are immediate.\\
In the sequel, we assume that $\zeta\neq 0$.\\
For the pressure $p$ we have then
\be
\label{5eq:37}
p_\zeta(t,z) =  q_\zeta(t) + \beta\int_0^z \theta_\zeta(t,z')\,dz'
\ee
so that existence and unicity for $u$, $v$ and $\theta$ give those of $p$ (up to the constant $c(t)$)  and(\ref{eq:16}).

\paragraph{Fourier-Laplace transform in time.}
\label{5sec:red2}
For $f(t,z) \in L^2(\mbb{R} ;H^\sigma_\zeta)$ with support in $\{t\geq 0\}$, we denote by $\hat{f}(\tau)$ its Fourier-Laplace transform:
\be
\label{5eq:38}
\hat{f}(\tau,z) = \int_0^{+\infty} e^{-it\tau}\, f(t,z)\,dt 
\ee
It is clear that  $\hat{f}$ is holomorphic in $\{\tau\in\mbb{C}, \Im (\tau)<0 \}$ and satisfies
\be
\label{5eq:39}
 \int_{-\infty}^{+\infty} \|\hat{f}(\tau)\|^2_{H^\sigma_\zeta}\,d\tau \,= \,C_0 \,\|f\|^2_{L^2(\mbb{R}_+;H^\sigma_\zeta)}
\ee
From  (\ref{5eq:25}), for a given $\zeta\neq 0$, for $u_\zeta, v_\zeta \in L^2(\mbb{R} ;H^\sigma_\zeta)$ we have
\be
\label{5eq:40}
p_\zeta \in L^2(\mbb{R}_+;H^\sigma_\zeta) + \dt L^2(\mbb{R}_+;H^{\sigma+2}_\zeta)
\ee
where $\dt L^2(\mbb{R}_+;H^{\sigma+2}_\zeta) = \{ q, \exists \tld{q} \in L^2(\mbb{R}_+;H^{\sigma+2}_\zeta), q = \dt \tld{q} \}$. Thus the Fourier-Laplace transform of $p_\zeta$ is well-defined.

\paragraph{Introduction of the parameters. } 
Let $\mbb{S}$ be the subset of $\mbb{C}$ defined by:
\be
\label{5eq:44}
\mbb{S} = \{-\delta_2 -\mu_1+i\mu_2,\, \textrm{ with } (\mu_1,\mu_2) \in \mbb{R}^2 \textrm{ and } |\mu_2|\geq \frac{1}{\delta_1} \mu_1 \}
\ee
with $\delta_1>0$ small enough  that $\mbb{S}\cap \mbb{V}_P = \emptyset$ (which is possible from lemma \ref{lem:3}) and
\be
\label{5eq:63}
\delta_2 < \delta_3 = \min (\frac{\nu\pi^2}{2a^2},\frac\nu2)
\ee
For $\zeta^2 \in \mbb{Z}^2\setminus 0$ and $\lambda\in\mbb{S}$, we set
\be
\label{5eq:41}
\lambda=i\tau, \quad \omega^2=\lambda + \nu\zeta^2
\ee
so that
\be
\label{5eq:70}
\left.\ba{l}\lambda\in\mbb{S}\\\zeta\in\mbb{Z}^2\setminus 0\ea\right\}
\quad\Rightarrow\quad
\left\{\ba{l} w^2\neq 0\\\lambda+\delta_3\neq 0\\\nu\zeta^2-\delta_3 >0 \ea\right.
\ee
Keeping the same notations for functions and their Fourier-Laplace transform, we obtain that  (\ref{5eq:25}) is equivalent to the following equations, with parameter  $(\zeta^2,\lambda) \in \mbb{Z}^2\setminus 0 \times\mbb{S}$:
\be
\label{5eq:62}
\ba{l}
\baa{l}
(\omega^2 -\nu \dzz) u_\zeta -\alpha v_\zeta +i\xi p_\zeta = f_{1,\zeta}\\
(\omega^2 -\nu \dzz) v_\zeta +\alpha u_\zeta +i\eta p_\zeta = f_{2,\zeta}\\
\dz p_\zeta - \beta \theta_\zeta = 0\\
(\omega^2 -\nu \dzz)\theta_\zeta+\gamma w_\zeta =f_{3,\zeta} 
\eaa \smallbreak\\
\textrm{ with }w_\zeta(z)=-\int_0^z (i\xi u_\zeta +i\eta v_\zeta)\\
\textrm{ and } w_\zeta(a)=0,\,\, X_\zeta|_{z=0,z=a}=0
\ea
\ee
For $\zeta$ and $\lambda$ given $\mbb{Z}^2\setminus 0$ and $\mbb{S}$, (\ref{5eq:62}) is a differential system as a function of $z\in(0,a)$ with data $F_\zeta\in H^\sigma(0,a)$ and unknown  $X_\zeta \in H^\sigma(0,a)$, $X_\zeta|_{z=0,a}=0$. Therefore, unicity is clear. Indeed, the third equation of  (\ref{5eq:62}) gives $p_\zeta\in H^{\sigma+3}(0,a)$, so that if $F=0$ we obtain $u_\zeta, v_\zeta \in H^{\sigma+5}(0,a)$, thus $w\in H^{\sigma+6}(0,a)$, then $\theta \in H^{\sigma+8}(0,a)$. In particular we have $X_\zeta \in H^1_0(0,a)$, and the spectral result $\mbb{S}\cap \mbb{V}_P = \emptyset$ gives $\hat{X}_\zeta=0$ for all $(\lambda,\zeta)\in\mbb{Z}^2\setminus 0\times\mbb{S}$. Then $X=0$ and $p=0$ up to a function of time.\\
It remains to prove existence and estimates for the solutions of  (\ref{5eq:62}), which will give (\ref{5eq:29}) and (\ref{5eq:219})  thanks to (\ref{5eq:39}).

\subsection{Preliminary results for the parameters} \label{thm1:2}
\label{5sec:param}
We set
\be
\label{5eq:51}
\la\zeta\ra = 1+ |\zeta|, \quad  \la\omega\ra^2=|\lambda|+\la\zeta\ra^2 
\ee
The following lemma will be useful:
\begin{lem} \label{5lem:5}
There exists a constant $C$ such that, for all $\lambda\in\mbb{S}$, $\zeta\in\mbb{Z}^2\setminus 0$ and $k\in\mbb{N}^*$:
\be
\label{5eq:201}
|\omega^2|\, \geq\, C\,\la\omega\ra^2\, \geq\, C\, (1+|\lambda|+\zeta^2)
\ee
and
\be
\label{5eq:202}
|\omega^2+\frac{\nu k^2\pi^2}{a^2}| \geq C\,(\la\omega\ra^2 + k^2).
\ee
\end{lem}
This lemma is easily obtained from the straightforward following lemma:
\begin{lem} \label{5lem:6}
Let $C_1$ and $C_2$ be two closed cones of $\mbb{R}^n$. We assume the distance between $C_1$ and $C_2$ to be non-zero, ie there exists a constant $d>0$ such that
\be
\label{5eq:203}
\forall x \in C_1,\forall y \in C_2,\quad \|x\|=\|y\|=1 \Rightarrow \|x-y\| \geq d
\ee
Then there exists a constant $C>0$ such that
\be
\label{5eq:204}
\forall x\in C_1, \forall y\in C_2,\quad \|x-y\| \geq C \,(\|x\|+\|y\|)
\ee
\end{lem}
\paragraph{The function $M_\sigma(\lambda,\zeta)$.}
For $(\lambda,\zeta)\in\mbb{S}\times\mbb{Z}^2\backslash 0$ we set:
\be
\label{5eq:78}
M_\sigma(\lambda,\zeta) = \Big(\, \sum_{k\in\mbb{N}^*} \frac1{k^2 (k^4+\la\omega\ra^4)(k^2+\la\zeta\ra^2)^\sigma}\, \Big)^{1/2}
\ee
which is well-defined for $\sigma>-\frac52$.\\
We shall use the following notation
\be
\label{5eq:91}
A(\lambda,\zeta) \sim B(\lambda,\zeta)\quad\Leftrightarrow\quad
\baa{l} \exists C_1,C_2>0, \, \forall (\lambda,\zeta)\in\mbb{S}\times\mbb{Z}^2\backslash 0,\medbreak\\
 \, C_1\, B(\lambda,\zeta)\,\,\leq\,\, A(\lambda,\zeta)\,\,\leq\,\, C_2\, B(\lambda,\zeta) 
\eaa
\ee
\begin{lem}\label{5prop:1}
For $\sigma\in]-\frac52,\frac32[$, $\sigma\neq -\frac12$, the following hold true:
\be
\label{5eq:88}
\ba{lccl}
\sigma>-\frac12 \,:\quad& M_\sigma& \sim& \dfrac{\la\zeta\ra^{-\sigma}}{\la\omega\ra^2}\medbreak\\
\sigma<-\frac12 \,:& M_\sigma& \sim&  \dfrac{\la\zeta\ra^{-\sigma}}{\la\omega\ra^2} \mbbm{1}_{\la\zeta\ra\geq\la\omega\ra^\kappa} +  \dfrac{\la\omega\ra^{-\sigma}}{\la\omega\ra^\frac52} \mbbm{1}_{\la\zeta\ra\leq\la\omega\ra^\kappa}\,, \quad \textrm{ with }\kappa=\frac{2\sigma+1}{2\sigma}
\ea
\ee
\end{lem}
\bproof First we can write: 
\be
\label{5eq:89}
\ba{rcl}
M_\sigma^2=\sum \frac{(k^2+\la\zeta\ra^2)^{-\sigma }}{ k^2\,(\la\omega\ra^4+ k^4)} &\sim& 
 \big[ \int_1^{\la\zeta\ra}  \frac{(x^2+\la\zeta\ra^2)^{-\sigma }}{ x^2\,(\la\omega\ra^4+x^4)}\,dx  
+ \int_{\la\zeta\ra}^{\la\omega\ra}\frac{(x^2+\la\zeta\ra^2)^{-\sigma }}{ x^2\,(\la\omega\ra^4+x^4)}\,dx  \smallbreak\\
&&\qquad+ \int_{\la\omega\ra}^{\infty} \frac{(x^2+\la\zeta\ra^2)^{-\sigma }}{ x^2\,(\la\omega\ra^4+x^4)}\,dx  \big]
\ea
\ee
We have immediately:
\be
\label{5eq:90}
\ba{lrcl}
I_1\,=&\int_1^{\la\zeta\ra}  \frac{(x^2+\la\zeta\ra^2)^{-\sigma }}{ x^2\,(\la\omega\ra^4+x^4)}\,dx  &\sim&  \frac{\la\zeta\ra^{-2\sigma}}{\la\omega\ra^4} \int_1^{\la\zeta\ra} \frac{dx}{x^2}  \smallbreak\\
&&\sim& \frac{\la\zeta\ra^{-2\sigma} }{\la\omega\ra^4}\bigbreak\\
I_2\,=& \int_{\la\zeta\ra}^{\la\omega\ra}\frac{(x^2+\la\zeta\ra^2)^{-\sigma }}{ x^2\,(\la\omega\ra^4+x^4)}\,dx &\sim&  \int_{\la\zeta\ra}^{\la\omega\ra}\frac{x^{-2\sigma }}{ x^2\,(\la\omega\ra^4+x^4)}\,dx \smallbreak\\
&&\sim& \frac{1}{\la\omega\ra^4}\int_{\la\zeta\ra}^{\la\omega\ra} \frac{dx}{x^{2\sigma+2}}\bigbreak\\
I_3\,=& \int_{\la\omega\ra}^{\infty} \frac{(x^2+\la\zeta\ra^2)^{-\sigma }}{ x^2\,(\la\omega\ra^4+x^4)}\,dx &\sim&\int_{\la\omega\ra}^{\infty}  \frac{x^{-2\sigma }}{ x^2\,(\la\omega\ra^4+x^4)}\,dx \smallbreak\\
&&\sim& \frac{\la\omega\ra^{-2\sigma}}{\la\omega\ra^5}\int_{1}^{\infty}  \frac{u^{-2\sigma }}{ u^2\,(1+ u^4)}\,du \smallbreak\\
&&\sim& \frac{\la\omega\ra^{-2\sigma}}{\la\omega\ra^5}
\ea
\ee
Then there are two cases:\\
$\star$ Case $\sigma>-\frac12$, with $\sigma<0$: then $\frac{\la\omega\ra^{-2\sigma}}{\la\omega\ra}$ is bounded and $\la\zeta\ra^{-2\sigma}\geq 1$ thus
\be
\label{5eq:92}
 \frac{\la\omega\ra^{-2\sigma}}{\la\omega\ra^5}\,\leq\,C \frac{1}{\la\omega\ra^4}\,\leq\,C \frac{\la\zeta\ra^{-2\sigma} }{\la\omega\ra^4}\quad\Rightarrow\quad I_3\,\leq\,C I_1
\ee
If $\sigma>-\frac12$, with $\sigma>0$, then $\la\omega\ra^{-2\sigma}\leq\la\zeta\ra^{-2\sigma}$ thus $I_3\,\leq\,C I_1$.\\
Moreover, if $\sigma>-\frac12$ we have $2\sigma+2>1$ then
\be
\label{5eq:93}
\int_{\la\zeta\ra}^{\la\omega\ra} \frac{dx}{x^{2\sigma+2}} \,\leq\, C \la\zeta\ra^{-2\sigma-1}\quad\Rightarrow\quad I_2\,\leq\,C I_1
\ee
It follows that  (\ref{5eq:88}) holds true for $\sigma>-\frac12$.\\
$\star$ Case $\sigma<-\frac12$: we have
\be
\label{5eq:106}
\int_{\la\zeta\ra}^{\la\omega\ra} \frac{dx}{x^{2\sigma+2}} \,\leq\, C \la\omega\ra^{-2\sigma-1}\quad\Rightarrow\quad I_2\,\leq\,C I_3
\ee
In order to compare $I_1$ and $I_3$ we introduce the following critical exponent: 
\be
\label{5eq:110}
\kappa \, =\, \frac{2\sigma+1}{2\sigma}\quad\Leftrightarrow\quad  -2\sigma\kappa\,=\,-2\sigma-1
\ee
If $\sigma\in]-\frac52,-\frac12[$, then  $\kappa \in ]0,\frac45[$.  We have again two cases:\\
$\ast$ If $(\la\zeta\ra,\la\omega\ra) \in \{\la\zeta\ra\leq\la\omega\ra^{\kappa} \}$, then
\be
\label{5eq:114}
\la\zeta\ra^{-2\sigma}\, \leq\, \la\omega\ra^{-2\sigma\kappa}\, =\, \la\omega\ra^{-2\sigma-1}\quad \Rightarrow \quad I_1 \,\leq\,C\, I_3
\ee
$\ast$ If $(\la\zeta\ra,\la\omega\ra) \in \{\la\zeta\ra\geq\la\omega\ra^{\kappa} \}$, then
\be
\label{5eq:125}
\la\zeta\ra^{-2\sigma}\, \geq\, \la\omega\ra^{-2\sigma\kappa}\, =\, \la\omega\ra^{-2\sigma-1}\quad \Rightarrow \quad I_3 \,\leq\,C\, I_1
\ee
And (\ref{5eq:88}) follows, for $\sigma<-\frac12$.\\
\eproof\\
We have also
\begin{cor}\label{5cor:2}
If $\sigma\in]-\frac12,\frac12[$ then
\be
\label{5eq:268}
\la\omega\ra^2 M_\sigma M_{-\sigma} \quad\sim\quad    \la\omega\ra^{-2}
\ee
If $\sigma\in]-\frac32,-\frac12[$ then
\be
\label{5eq:269}
\la\omega\ra^2 M_\sigma M_{-\sigma} \quad\leq\quad \la\omega\ra^{-1}
\ee
and
\be
\label{5eq:285}
\la\omega\ra^2 M_\sigma M_{-\sigma-2} \quad \leq \quad C\, \mbbm{1}_{\la\zeta\ra\geq\la\omega\ra^\kappa} + \la\omega\ra^{-1}\mbbm{1}_{\la\zeta\ra\leq\la\omega\ra^\kappa}
\ee
\end{cor}
\bproof If $\sigma\in]-\frac12,\frac12[$, (\ref{5eq:268}) is immediately obtained from lemma \ref{5prop:1}. If $\sigma\in]-\frac32,-\frac12[$ then
\be
\label{5eq:284}
\ba{rcl}
\la\omega\ra^2 M_\sigma M_{-\sigma} &\sim& \la\omega\ra^2 \big(\frac{\la\zeta\ra^{-\sigma}}{\la\omega\ra^2} \mbbm{1}_{\la\zeta\ra\geq\la\omega\ra^\kappa} + \frac{\la\omega\ra^{-\sigma}}{\la\omega\ra^{5/2}} \mbbm{1}_{\la\zeta\ra\leq\la\omega\ra^\kappa} \big) \frac{\la\zeta\ra^{\sigma}}{\la\omega\ra^2} \smallbreak\\
&\sim& \la\omega\ra^{-2} \mbbm{1}_{\la\zeta\ra\geq\la\omega\ra^\kappa} +\la\zeta\ra^{\sigma}\la\omega\ra^{-\sigma-5/2} \mbbm{1}_{\la\zeta\ra\leq\la\omega\ra^\kappa} \smallbreak\\
&\leq& \la\omega\ra^{-2} \mbbm{1}_{\la\zeta\ra\geq\la\omega\ra^\kappa} +\la\omega\ra^{-\sigma-5/2} \mbbm{1}_{\la\zeta\ra\leq\la\omega\ra^\kappa} \smallbreak\\
&\leq& \la\omega\ra^{-1}
\ea
\ee
because $\la\zeta\ra^{\sigma}<1$ as $\sigma<0$ and $\la\omega\ra^{-\sigma-5/2} < \la\omega\ra^{-1}$ as $-\sigma-5/2<-1$.\\
Finally we have
\be
\label{5eq:283}
\ba{rcl}
\la\omega\ra^2 M_\sigma M_{-\sigma-2} &\sim& \la\omega\ra^2 (\la\zeta\ra^{-\sigma}\la\omega\ra^{-2}\mbbm{1}_{\la\zeta\ra\geq\la\omega\ra^\kappa} + \la\omega\ra^{-\sigma-5/2}\mbbm{1}_{\la\zeta\ra\leq\la\omega\ra^\kappa}   ) \smallbreak\\
&&\qquad\qquad ( \la\zeta\ra^{\sigma+2}\la\omega\ra^{-2}\mbbm{1}_{\la\zeta\ra\geq\la\omega\ra^\kappa} + \la\omega\ra^{\sigma-1/2}\mbbm{1}_{\la\zeta\ra\leq\la\omega\ra^\kappa}    )\medbreak\\
&\sim & \la\zeta\ra^{2}\la\omega\ra^{-2}\mbbm{1}_{\la\zeta\ra\geq\la\omega\ra^\kappa} + \la\omega\ra^{-1}\mbbm{1}_{\la\zeta\ra\leq\la\omega\ra^\kappa}  \medbreak\\
&\leq&C\, \mbbm{1}_{\la\zeta\ra\geq\la\omega\ra^\kappa} + \la\omega\ra^{-1}\mbbm{1}_{\la\zeta\ra\leq\la\omega\ra^\kappa}
\ea
\ee
\eproof

\subsection{The uncoupled system}\label{thm1:3}
\label{5sec:dec}
Having disposed of these preliminary steps, we can now study the uncoupled system (\ref{5eq:62}) ie with $\alpha=\beta=\gamma=0$:
\be
\label{5eq:126}
\ba{l}
\baa{l}
(\omega^2 -\nu \dzz) u_\zeta +i\xi p_\zeta = f_{1,\zeta}\\
(\omega^2 -\nu \dzz) v_\zeta +i\eta p_\zeta = f_{2,\zeta}\\
\dz p_\zeta  = 0\\
(\omega^2 -\nu \dzz)\theta_\zeta =f_{3,\zeta} 
\eaa \smallbreak\\
\textrm{ with }\int_0^a (\xi u_\zeta +\eta v_\zeta)=0\textrm{ and }  X_\zeta|_{z=0,z=a}=0
\ea
\ee

\paragraph{Notations.} We denote for short:
\begin{itemize}
\item $p_{\zeta,0}=p_\zeta(z=0)$;
\item $\|.\|_{\sigma,\zeta} $ stands for $\|.\|_{H^{\sigma}_\zeta}$, $\|.\|_{(H^{\sigma}_\zeta)^2}$  and  $\|.\|_{(H^{\sigma}_\zeta)^3}$  ;
\item $(H^{\sigma+2}_\zeta)^2_\divg$ stands for $\{ (u,v)\in (H^{\sigma+2}_\zeta)^2,\int_0^a \xi u +\eta v=0\}$.\\
\end{itemize}
The next theorem is the core of the proof:
\begin{thm}\label{5prop:3}
Let $(\lambda,\zeta)\in\mbb{S}\times\mbb{Z}^2\backslash 0$. The operator
\be
\label{5eq:128}
\ba{cccc}
& (H^{\sigma+2}_\zeta)^2_\divg \times \mbb{C}\times H^{\sigma+2}_\zeta &\rightarrow& (H^{\sigma}_\zeta)^3\medbreak\\
\mc{L}_0:& (u,v,p_0,\theta)&\mapsto& \left[\ba{l} (\omega^2 -\nu \dzz) u +i\xi p_0\\ (\omega^2 -\nu \dzz) v +i\eta p_0\\(\omega^2 -\nu \dzz)\theta \ea\right] = \left[\ba{l} f_1\\f_2\\f_3 \ea\right] 
\ea
\ee
is continuous and bijective. Moreover $Y=(u,v)$ splits in $Y_1+Y_2$ and the following estimates hold:
\be
\label{5eq:129}
\ba{lrcl}
(a)&|\zeta p_0| &\leq& C\la\omega\ra^2 M_\sigma \|F\|_{\sigma,\zeta} \medbreak\\
(b)&M_{-\sigma-2}[M_{-\sigma}]^{-1}\|Y_1\|_{\sigma,\zeta} + \|Y_1\|_{\sigma+2,\zeta} &\leq& C \la\omega\ra^2 M_\sigma M_{-\sigma-2}\|F\|_{\sigma,\zeta} \medbreak\\
(c)&\la\omega\ra^2 \|Y_2\|_{\sigma,\zeta} + \|Y_2\|_{\sigma+2,\zeta} &\leq& C  \|F\|_{\sigma,\zeta}\medbreak\\
(d)& \la\omega\ra^2 \|\theta\|_{\sigma,\zeta} +  \|\theta\|_{\sigma+2,\zeta} &\leq& C \|f_{3}\|_{\sigma,\zeta}
\ea
\ee
with $(f_1,f_2)=F$.
\end{thm}
Before proving this theorem, we state an elementary but useful remark:
\begin{rmq}\label{5rmq:11}
The eigenvalues of the operator $\omega^2-\nu\dzz$ are eigenvalues of  $P$. Moreover, for all $\lambda\notin\mbb{V}_P$, the operator $\omega^2-\nu\dzz$ is continuous and bijective from  $H^{\sigma+2}_\zeta$ to $H^{\sigma}_\zeta$, for all $\sigma\in]-\frac32,\frac12[$.
\end{rmq}

\paragraph{Proof of theorem\ref{5prop:3}. } Continuity is easy, as is injectivity (because $\mbb{S}\cap\mbb{V}_P=\emptyset$). \\
We first examine the $\theta$ part. Surjectivity is clear. Then $\theta$ satisfies a heat equation for which we can write a maximal estimate, which is exactly (\ref{5eq:129},d): indeed let $f\in H^{\sigma}_\zeta$, with $f=\sum_{k\in\mbb{N}} f_k e_k(z)$. Let $g$ be the solution of  $(\omega^2 -\nu \dzz)g=f$ with homogeneous Dirichlet boundary conditions. Then $g\in H^{\sigma+2}(0,a)$, with $\sigma+2 \in ]\frac12 ,\frac52[$, and $g(0)=g(a)=0$, hence  $g\in H^{\sigma+2}_\zeta$ according to lemma \ref{5lem:4}. Thus $g=\sum_k g_k e_k$ and we have:
\be
\label{5eq:156}
(\omega^2+\frac{\nu k^2\pi^2}{a^2}) g_k = f_k
\ee
with $\omega^2+\frac{\nu k^2\pi^2}{a^2}\neq 0$.\\
Using lemma \ref{5lem:5} we obtain easily:
\be
\label{5eq:162}
\ba{rcl}
\la\omega\ra^4 \|g\|^2_{\sigma,\zeta} &=& \la\omega\ra^4 \dsum_k (1+\nu\zeta^2+\nu k^2)^\sigma \dfrac{|f_k|^2}{|\omega^2+\frac{\nu k^2\pi^2}{a^2}|^2}\smallbreak\\
&\leq& C\la\omega\ra^4 \dsum_k (1+\nu\zeta^2+\nu k^2)^\sigma \dfrac{|f_k|^2}{\la\omega\ra^4+k^4} \quad\leq\quad C \|f\|^2_{\sigma,\zeta}
\ea
\ee
and
\be
\label{5eq:244}
\ba{rcl}
\|g\|^2_{\sigma+2,\zeta} &\leq& \dsum_k (1+\nu\zeta^2+\nu k^2)^{\sigma} |f_k|^2\dfrac{(1+\nu\zeta^2+\nu k^2)^2}{|\omega^2+\frac{\nu k^2\pi^2}{a^2}|^2}\smallbreak\\
&\leq&C \dsum_k (1+\nu\zeta^2+\nu k^2)^{\sigma} |f_k|^2\dfrac{\la\zeta\ra^4+k^4}{\la\omega\ra^4+k^4}\quad \leq \quad C \|f\|^2_{\sigma,\zeta}
\ea
\ee
With $f=f_{3,\zeta}\in H^{\sigma}_\zeta$ and $\theta_\zeta=g$ we get (\ref{5eq:129},d).\\
We know examine $\mc{L}_0$ restricted to $(u,v,p_0)\in (H^{\sigma+2}_\zeta)^2_\divg \times \mbb{C}$. To establish surjectivity, we first evaluate (explicitly) the constant $p_0$, we then invert the operator  $(\omega^2 -\nu \dzz)$ (which is possible according to remark \ref{5rmq:11}). Thus it is sufficient to evaluate $p_0$ and prove estimates (\ref{5eq:129}).\\ 
With (\ref{5eq:128}) we obtain:
\be
\label{5eq:130}
u=(\omega^2 -\nu \dzz)^{-1}[f_1-i\xi p_0],\quad v=(\omega^2 -\nu \dzz)^{-1}[f_1-i\eta p_0]
\ee
then we evaluate $p_0$ thanks to $\int_0^ai \xi u +i\eta v=0$:
\be
\label{5eq:131}
\ba{rl}
&\int_0^a i\xi (\omega^2 -\nu \dzz)^{-1}[f_1-i\xi p_0] +i\eta(\omega^2 -\nu \dzz)^{-1}[f_1-i\eta p_0]\,dz =0\smallbreak\\
\Leftrightarrow & \int_0^a(\omega^2 -\nu \dzz)^{-1}[(\xi^2+\eta^2 )p_0]\,dz =- \int_0^a  (\omega^2 -\nu \dzz)^{-1}[i\xi f_1+i\eta f_2]  \,dz\smallbreak\\
\Leftrightarrow & \zeta^2 p_0 \int_0^a(\omega^2 -\nu \dzz)^{-1}[1]\,dz =- \int_0^a  (\omega^2 -\nu \dzz)^{-1}[i\xi f_1+i\eta f_2]  \,dz
\ea
\ee
We will prove further the following lemma:
\begin{lem}\label{5lem:7} There exist constants $C_1$ and $C_2$, independent of $\omega^2=\lambda+\nu\zeta$, such that, for all $(\lambda, \zeta)\in \mbb{S}\times \mbb{Z}\setminus 0$, 
\be \label{eq:25}
\frac{C_1}{\langle \omega \rangle^2} \,\leq\,  \dint_0^a (\lambda -\nu\Delta)^{-1}(1 )\, dz   \,\leq\,\frac{C_2}{\langle \omega \rangle^2}
\ee
\end{lem}
Thus $p_0$ is well-defined by the following formula:
\be
\label{5eq:155}
\zeta^2 p_0  =- \Big[   \dint_0^a (\lambda -\nu\Delta)^{-1}(1 )\, dz   \Big]^{-1} \int_0^a  (\omega^2 -\nu \dzz)^{-1}[i\xi f_1+i\eta f_2]  \,dz
\ee
Let us now estimate $\int_0^a  (\omega^2 -\nu \dzz)^{-1}[i\xi f_1+i\eta f_2]  \,dz$. Let $f=\sum_k f_k e_k \in H^\sigma_\zeta$. As previously, using lemma \ref{5lem:5}, we have
\be
\label{5eq:157}
\ba{rcl}
|\int_{0}^a (\omega^2 -\nu \dzz)^{-1} f| &=& |\int_0^a \sum \dfrac{f_k}{\omega^2+\frac{\nu k^2\pi^2}{a^2}} e_k(z)\, dz|\smallbreak\\
&\leq& C \sum\dfrac{|f_k|}{k\,(\la\omega\ra^2+ k^2)}\smallbreak\\
&\leq& C \|f\|_{\sigma,\zeta}  \Big[\sum \dfrac{(1+\nu k^2+\nu|\zeta|^2)^{-\sigma }}{ k^2\,(\la\omega\ra^2+ k^2)^2}\Big]^{1/2}\smallbreak\\
&\leq& C \|f\|_{\sigma,\zeta}  M_\sigma(\lambda,\zeta)
\ea
\ee
Using (\ref{eq:25}), (\ref{5eq:155}) and (\ref{5eq:157}) we obtain:
\be
\label{5eq:158}
\zeta^2 \, |p_0| \, \leq \, C \, \la\omega\ra^2  \, M_\sigma(\lambda,\zeta)  \,|\zeta| \, \|F\|_{\sigma,\zeta}
\ee
and  (\ref{5eq:129},a) and surjectivity follows.\\
To establish (\ref{5eq:129},b,c) we use (\ref{5eq:130}):
\be
\label{5eq:159}
u =   (\omega^2-\nu\dzz)^{-1}[f_1] -i\xi p_0 (\omega^2-\nu\dzz)^{-1}[1]= u_2 + u_1
\ee
Define
\be
\label{5eq:160}
\ba{rcl}
Y_1 &=& -i \left[ \ba{c} \xi\\ \eta \ea \right] p_0 (\omega^2-\nu\dzz)^{-1}[1]\medbreak\\
Y_2 &=& (\omega^2-\nu\dzz)^{-1}\left[ \ba{c} f_1\\ f_2 \ea \right]
\ea
\ee
Like $\theta_\zeta$, $Y_2$ satisfies the maximal estimate of the heat equation, namely:
\be
\label{5eq:161}
\la\omega\ra^2 \|Y_2\|_{\sigma,\zeta} +  \|Y_2\|_{\sigma+2,\zeta} \leq C \|F\|_{\sigma,\zeta}
\ee
We now turn to $Y_1$:
\be
\label{5eq:164}
\|Y_1\|_{\sigma+2,\zeta} = |\zeta| \,|p_0| \, \|(\omega^2-\nu\dzz)^{-1}[1] \|_{\sigma+2,\zeta}
\ee
As previously, the constant function $1$ stays in $H^{\sigma}_\zeta$ and $1 = \sum_k a_k e_k$, with $a_k = C_0 \frac1k$ for $k$ odd and $a_k=0$ otherwise, and we obtain:
\be
\label{5eq:168}
\|(\omega^2-\nu\dzz)^{-1}[1] \|_{\sigma+2,\zeta} 
~\leq~ C\Big[ \dsum_{k}  \dfrac{(\la\zeta\ra^2+ k^2)^{\sigma+2}}{k^2(\la\omega\ra^4+k^4)} \Big]^{1/2} ~=~ C M_{-\sigma-2}(\lambda,\zeta)
\ee
($M_{-\sigma-2}$ is well-defined for $\sigma \in ]-\frac32, \frac12[$)\\
Likewise we get:
\be
\label{5eq:199}
\|(\omega^2-\nu\dzz)^{-1}[1] \|_{\sigma,\zeta} \,\leq\, C M_{-\sigma}(\lambda,\zeta)
\ee
($M_{-\sigma}$ is well-defined for  $\sigma \in ]-\frac32, \frac12[$)\\
which completes the proof.\\
\eproof

\paragraph{Proof of lemma \ref{5lem:7}. } Straightforward calculation gives:
\be
\label{5eq:148}
\dint_0^a (\lambda -\nu\Delta)^{-1}(1 )\, dz =  \dfrac{a}{\omega^2} \Big[ 1 + 2\dfrac{1- \cosh(\frac{\omega a}{\sqrt{\nu}})}{\frac{\omega a}{\sqrt{\nu}}\sinh(\frac{\omega a}{\sqrt{\nu}})} \Big]
\ee
where $\omega$ is the square root of $\omega^2$ with positive real part. Let us define
\be
\label{5eq:147}
\mc{N}(\chi) =1 + 2 \dfrac{1-\cosh(\chi)}{\chi\sinh(\chi)} 
\ee
with $\chi=\frac{\omega a}{\sqrt{\nu}}$. We have $\omega^2=\lambda+\nu\zeta$, thus according to (\ref{5eq:44}), (\ref{5eq:63}) and (\ref{5eq:70}) 
$\omega$ and $\chi$ stays in $\mbb{B}$ with
\be
\label{5eq:154}
\mbb{B}= \{\mu_1+i\mu_2, \textrm{ with } (\mu_1,\mu_2) \in \mbb{R}_+\times\mbb{R} \textrm{ and } |\mu_2|\leq (1+\delta_5) \mu_1 \} \setminus \delta_4 B_1 
\ee 
where $\delta_4$ and $\delta_5$ are small and $B_1$ is the open unity disk. \\
The singularities of $\mc{N}$ are zeros of $\sinh$, so that  $\mc{N}$ is holomorphic in $\mbb{B}$. Moreover, $\mc{N}$ has no zero in $\mbb{B}$, because  $\mc{L}_0$ is injective for $\omega\in\mbb{B}$. We have besides
\be
\label{5eq:151}
\lim_{|\chi|\rightarrow +\infty}\mc{N}(\chi)=1
\ee
Hence there exist  constants  $C_1$ and $C_2$, independent of $\omega$, such that for all  $\chi\in\mbb{B}$
\be
\label{5eq:149}
\frac{C_1}{1+|\chi|^2} \,\leq\, \frac{|\mc{N}(\chi)|}{|\chi|^2} \,\leq\,\frac{C_2}{1+|\chi|^2}
\ee
Estimate (\ref{eq:25}) easily follows from (\ref{5eq:149}) and lemma \ref{5lem:5}.\\
\eproof

\subsection{The coupled system}\label{thm1:4}
\label{5sec:coup}
We now turn to the case $\alpha,\beta,\gamma\neq 0$ and system (\ref{5eq:62}). Let us define the corresponding perturbation operator $\mc{L}_1$:
\be
\label{5eq:173}
\ba{cccc}
& (H^{\sigma+2}_\zeta)^2_\divg \times \mbb{C}\times H^{\sigma+2}_\zeta &\rightarrow& (H^{\sigma}_\zeta)^3\medbreak\\
\mc{L}_1:& (u,v,p_0,\theta)&\mapsto& \left[\ba{l}  -\alpha v +i\xi \beta \int_0^z \theta \smallbreak\\ \alpha u+i\eta \beta \int_0^z \theta\smallbreak\\-\gamma \int_0^z( i\xi u+i\eta v) \ea\right] 
\ea
\ee
such that $\mc{L}=\mc{L}_0+\mc{L}_1$, where $\mc{L}_0$ is given by (\ref{5eq:128}), corresponds to system (\ref{5eq:62}). We have the
\begin{thm}\label{5prop:8}
Let $(\lambda,\zeta)\in\mbb{S}\times\mbb{Z}^2\backslash 0$. The following operator
\be
\label{5eq:175}
\ba{cccc}
& (H^{\sigma+2}_\zeta)^2_\divg \times \mbb{C}\times H^{\sigma+2}_\zeta &\rightarrow& (H^{\sigma}_\zeta)^3\medbreak\\
\mc{L}:& (u,v,p_0,\theta)&\mapsto& (\mc{L}_0+\mc{L}_1) (u,v,p_0,\theta)= F
\ea
\ee
is continuous and bijective. Moreover $Y=(u,v)$ splits in $Y_1+Y_2$ and the following estimates hold true:
\be
\label{5eq:174}
\ba{lrcl}
(a)&|\zeta p_0| &\leq& C\la\omega\ra^2 M_\sigma \|F\|_{\sigma,\zeta} \medbreak\\
(b)&M_{-\sigma-2}[M_{-\sigma}]^{-1}\|Y_1\|_{\sigma,\zeta} + \|Y_1\|_{\sigma+2,\zeta} &\leq& C \la\omega\ra^2 M_\sigma M_{-\sigma-2}\|F\|_{\sigma,\zeta} \medbreak\\
(c)&\la\omega\ra^2 \|Y_2\|_{\sigma,\zeta} + \|Y_2\|_{\sigma+2,\zeta} &\leq& C  \|F\|_{\sigma,\zeta}\medbreak\\
(d)& \la\omega\ra^2 \|\theta\|_{\sigma,\zeta} +  \|\theta\|_{\sigma+2,\zeta} &\leq& C \|F\|_{\sigma,\zeta}
\ea
\ee
\end{thm}
\bproof Let $(\lambda,\zeta)\in\mbb{S}\times\mbb{Z}^2\backslash 0$. To prove that $\mc{L}$ is an isomorphism, we first state that the image of  $\mc{L}_1$ is included in a compact subspace of  $(H^{\sigma}_\zeta)^3$. Indeed, if $(u,v,p_0,\theta)  \in (H^{\sigma+2}_\zeta)^2_\divg \times \mbb{C}\times H^{\sigma+2}_\zeta $, then 
\be
\label{5eq:176}
\ba{rcl}
-\alpha v +i\xi \beta \int_0^z \theta &\in& H^{\sigma+2}(0,a)\smallbreak\\
 \alpha u+i\eta \beta \int_0^z \theta &\in& H^{\sigma+2}(0,a)\smallbreak\\
-\gamma \int_0^z (i\xi u+i\eta v)  &\in& H^{\sigma+3}(0,a)
\ea
\ee
and $ (H^{\sigma+2}(0,a))^2\times  H^{\sigma+3}(0,a)$ is a compact subspace of $(H^{\sigma}_\zeta)^3$. From Fredholm theory,  $\mc{L}=\mc{L}_0+\mc{L}_1$ is an isomorphism if and only if its kernel is trivial. Let  $(u,v,p_0,\theta)$ be in the kernel of $\mc{L}$, ie
\be
\label{5eq:197}
\ba{l}
\baa{l}
(\omega^2 -\nu \dzz) u  = \alpha v -i\xi p\\
(\omega^2 -\nu \dzz) v =-\alpha u -i\eta p\\
\dz p - \beta \theta = 0\\
(\omega^2 -\nu \dzz)\theta=-\gamma w
\eaa \smallbreak\\
\textrm{ with }w(z)=-\int_0^z (i\xi u +i\eta v)\\
\textrm{ and } w(a)=0,\,\, X|_{z=0,z=a}=0
\ea
\ee
Thus $u$, $v$ and $\theta$ are smooth and   $0 \notin \mbb{V}_P$ gives $(u,v,\theta)=0$ and then $p_0=0$.\\
Let us turn now to estimates (\ref{5eq:174}). We write
\be
\label{5eq:198}
\mc{L}(u,v,p_0,\theta)=F\quad\Leftrightarrow\quad \mc{L}_0(u,v,p_0,\theta)=F-\mc{L}_1(u,v,p_0,\theta)
\ee
and we use theorem \ref{5prop:3}. We get:
\be
\label{5eq:245}
\ba{lrcl}
(a)&|\zeta p_0| &\leq& C\la\omega\ra^2 M_\sigma \big[\|F\|_{\sigma,\zeta} +\|Y\|_{\sigma,\zeta} +|\zeta|\|\int_0^z \theta\, dz \|_{\sigma,\zeta} \big]  \medbreak\\
(b)&\frac{M_{-\sigma-2}}{M_{-\sigma}}\|Y_1\|_{\sigma,\zeta} + \|Y_1\|_{\sigma+2,\zeta} &\leq& C \la\omega\ra^2 M_\sigma M_{-\sigma-2} \big[\|F\|_{\sigma,\zeta} +\|Y\|_{\sigma,\zeta} +|\zeta|\|\int_0^z \theta\, dz \|_{\sigma,\zeta} \big]\medbreak\\
(c)&\la\omega\ra^2 \|Y_2\|_{\sigma,\zeta} + \|Y_2\|_{\sigma+2,\zeta} &\leq& C \big[\|F\|_{\sigma,\zeta} +\|Y\|_{\sigma,\zeta} +|\zeta|\|\int_0^z \theta\, dz \|_{\sigma,\zeta} \big] \medbreak\\
(d)& \la\omega\ra^2 \|\theta\|_{\sigma,\zeta} +  \|\theta\|_{\sigma+2,\zeta} &\leq& C \big[\|F\|_{\sigma,\zeta} + |\zeta|\|\int_0^z Y\, dz \|_{\sigma,\zeta} \big]
\ea
\ee
We will prove further the following lemma:
\begin{lem}\label{5lem:8}
Let $s\in ]-\frac32,\frac12[$, $s\neq -\frac12$ and $\varphi\in H^s_\zeta\cap H^{s+2}_\zeta$. The function $\phi(z) = \int_0^z \varphi(z')\,dz'$ stays in $H^s_\zeta$.\\
Besides, if $s \in]-\frac12,\frac12[$, we get
\be
\label{5eq:246}
\|\phi\|_{s,\zeta} \,\leq \, C  \|\varphi\|_{s,\zeta}
\ee
where the constant $C$ is  independent of $\zeta$ and $\omega$.\\
If $s\in]-\frac32,-\frac12[$, we get
\be
\label{5eq:262}
\ba{rcl}
\|\phi\|_{s,\zeta} &\leq & C_2(\omega,\zeta) \, \big( \la\omega\ra^2 \|\varphi\|_{s,\zeta} + \|\varphi\|_{s+2,\zeta} \big)\medbreak\\
\|\phi\|_{s,\zeta} &\leq & C_1(\omega,\zeta) \, \big( \frac{M_{-s-2}}{M_{-s}} \|\varphi\|_{s,\zeta} + \|\varphi\|_{s+2,\zeta} \big)\ea
\ee
with:
\be
\label{5eq:263}
\ba{llcl}
(i)\qquad&  \la\zeta\ra  \, C_2(\omega,\zeta)  &\stackrel{\la\omega\ra \rightarrow \infty}{\longrightarrow}&0\medbreak\\
(ii)&\la\zeta\ra\, C_1(\omega,\zeta)\,\mbbm{1}_{\la\zeta\ra\geq\la\omega\ra^\kappa}  &\stackrel{\la\omega\ra \rightarrow \infty}{\longrightarrow}&0 \medbreak\\
(iii)& C_1(\omega,\zeta) \,\mbbm{1}_{\la\zeta\ra\leq\la\omega\ra^\kappa} &\stackrel{\la\omega\ra \rightarrow \infty}{\longrightarrow}&0 \ea
\ee
\end{lem}
Assume first that $\sigma \in]-\frac12,\frac12[$. Lemma \ref{5lem:8} gives:
\be
\label{5eq:208}
\ba{lrcl}
(a)&|\zeta p_0| &\leq& C\la\omega\ra^2 M_\sigma \big[\|F\|_{\sigma,\zeta} +\|Y\|_{\sigma,\zeta} +|\zeta|\|\theta\|_{\sigma,\zeta} \big]  \medbreak\\
(b)&\frac{M_{-\sigma-2}}{M_{-\sigma}}\|Y_1\|_{\sigma,\zeta} + \|Y_1\|_{\sigma+2,\zeta} &\leq& C \la\omega\ra^2 M_\sigma M_{-\sigma-2} \big[\|F\|_{\sigma,\zeta} +\|Y\|_{\sigma,\zeta} +|\zeta|\|\theta\|_{\sigma,\zeta} \big]\medbreak\\
(c)&\la\omega\ra^2 \|Y_2\|_{\sigma,\zeta} + \|Y_2\|_{\sigma+2,\zeta} &\leq& C \big[\|F\|_{\sigma,\zeta} +\|Y\|_{\sigma,\zeta} +|\zeta|\|\theta\|_{\sigma,\zeta} \big] \medbreak\\
(d)& \la\omega\ra^2 \|\theta\|_{\sigma,\zeta} +  \|\theta\|_{\sigma+2,\zeta} &\leq& C \big[\|F\|_{\sigma,\zeta} + |\zeta|\|Y\|_{\sigma,\zeta} \big]
\ea
\ee
Let us absorb the perturbative terms of the right-hand side. When $\la\omega\ra$ is bounded, estimates (\ref{5eq:174}) are true because $\mc{L}$ is an isomorphism. We shall then assume that $\la\omega\ra$ is large enough. Absorbing $\|Y_2\|_{\sigma,\zeta}$ in (\ref{5eq:208}$,c$) is easy and one gets: 
\be
\label{5eq:224}
\|Y_2\|_{\sigma,\zeta} \,\leq\, C \frac{1}{\la\omega\ra^2 } \big[\|F\|_{\sigma,\zeta} +\|Y_1\|_{\sigma,\zeta} +|\zeta|\|\theta\|_{\sigma,\zeta} \big]
\ee
so that (\ref{5eq:208}) becomes:
\be
\label{5eq:209}
\ba{lrcl}
(a)&|\zeta p_0| &\leq& C\la\omega\ra^2 M_\sigma \big[\|F\|_{\sigma,\zeta} +\|Y_1\|_{\sigma,\zeta} +|\zeta|\|\theta\|_{\sigma,\zeta} \big]  \medbreak\\
(b)&\frac{M_{-\sigma-2}}{M_{-\sigma}}\|Y_1\|_{\sigma,\zeta} + \|Y_1\|_{\sigma+2,\zeta} &\leq& C \la\omega\ra^2 M_\sigma M_{-\sigma-2} \big[\|F\|_{\sigma,\zeta} +\|Y_1\|_{\sigma,\zeta} +|\zeta|\|\theta\|_{\sigma,\zeta} \big]\medbreak\\
(c)&\la\omega\ra^2 \|Y_2\|_{\sigma,\zeta} + \|Y_2\|_{\sigma+2,\zeta} &\leq& C \big[\|F\|_{\sigma,\zeta} +\|Y_1\|_{\sigma,\zeta} +|\zeta|\|\theta\|_{\sigma,\zeta} \big] \medbreak\\
(d)& \la\omega\ra^2 \|\theta\|_{\sigma,\zeta} +  \|\theta\|_{\sigma+2,\zeta} &\leq& C \big[\|F\|_{\sigma,\zeta} + |\zeta|\|Y_1\|_{\sigma,\zeta}+\frac{|\zeta|^2}{\la\omega\ra^2}\|\theta\|_{\sigma,\zeta}  \big]
\ea
\ee
According to corollary \ref{5cor:2},  $\la\omega\ra^2 M_\sigma M_{-\sigma}\sim\la\omega\ra^{-2}$, therefore we can absorb $\|Y_1\|_{\sigma,\zeta}$ in (\ref{5eq:209}$,b$) to get:
\be
\label{5eq:210}
\|Y_1\|_{\sigma,\zeta} \,\leq\, C   \la\omega\ra^{-2}  \big[\|F\|_{\sigma,\zeta} +|\zeta|\|\theta\|_{\sigma,\zeta} \big]
\ee
so that (\ref{5eq:209}$,a,c$) gives:
\be
\label{5eq:211}
\ba{lrcl}
(a)&|\zeta p_0| &\leq& C\la\omega\ra^2 M_\sigma \big[\|F\|_{\sigma,\zeta}  +|\zeta|\|\theta\|_{\sigma,\zeta} \big]  \medbreak\\
(b)&\frac{M_{-\sigma-2}}{M_{-\sigma}}\|Y_1\|_{\sigma,\zeta} + \|Y_1\|_{\sigma+2,\zeta} &\leq& C \la\omega\ra^2 M_\sigma M_{-\sigma-2} \big[\|F\|_{\sigma,\zeta}  +|\zeta|\|\theta\|_{\sigma,\zeta} \big]\medbreak\\
(c)&\la\omega\ra^2 \|Y_2\|_{\sigma,\zeta} + \|Y_2\|_{\sigma+2,\zeta} &\leq& C \big[\|F\|_{\sigma,\zeta} +|\zeta|\|\theta\|_{\sigma,\zeta} \big] \medbreak\\
(d)& \la\omega\ra^2 \|\theta\|_{\sigma,\zeta} +  \|\theta\|_{\sigma+2,\zeta} &\leq& C \big[\|F\|_{\sigma,\zeta} +\|\theta\|_{\sigma,\zeta}  \big]
\ea
\ee
Absorbing $\|\theta\|_{\sigma,\zeta}$ is then easy and one gets:
\be
\label{5eq:212}
\|\theta\|_{\sigma,\zeta} \leq C \frac1{\la\omega\ra^2 }\|F\|_{\sigma,\zeta} 
\ee
so that (\ref{5eq:174}) is established and the proof is complete for $\sigma\in]-\frac12,\frac12[$.\\
Let us now turn to the case  $\sigma\in]-\frac32,-\frac12[$. Estimates (\ref{5eq:245})  and lemma  \ref{5lem:8} with $C_1$ for $Y_1$ and $C_2$ for $\theta$ and $Y_2$ give:
\be
\label{5eq:264}
\ba{lrcl}
(a)&|\zeta p_0| &\leq& C\la\omega\ra^2 M_\sigma \big[\|F\|_{\sigma,\zeta} +\|Y_1\|_{\sigma,\zeta} +\|Y_2\|_{\sigma,\zeta} +|\zeta| C_2 I_\theta \big]  \medbreak\\
(b)&I_1 &\leq& C \la\omega\ra^2 M_\sigma M_{-\sigma-2} \big[\|F\|_{\sigma,\zeta} +\|Y_1\|_{\sigma,\zeta}+\|Y_2\|_{\sigma,\zeta} +|\zeta| C_2 I_\theta  \big]\medbreak\\
(c)&I_2 &\leq& C \big[\|F\|_{\sigma,\zeta} +\|Y_1\|_{\sigma,\zeta}+\|Y_2\|_{\sigma,\zeta} +|\zeta| C_2 I_\theta  \big] \medbreak\\
(d)& I_\theta &\leq& C \big[\|F\|_{\sigma,\zeta} + |\zeta| C_2 I_2 + |\zeta| C_1 I_1 \big]
\ea
\ee
where we have denoted
\be
\label{5eq:265}
\ba{rcl}
I_\theta&=& \la\omega\ra^2 \|\theta\|_{\sigma,\zeta} +  \|\theta\|_{\sigma+2,\zeta}  \smallbreak\\
I_2&=& \la\omega\ra^2 \|Y_2\|_{\sigma,\zeta} +  \|Y_2\|_{\sigma+2,\zeta}  \smallbreak\\
I_1 &=&\frac{M_{-\sigma-2}}{M_{-\sigma}} \|Y_1\|_{\sigma,\zeta} +\|Y_1\|_{\sigma+2,\zeta} 
\ea
\ee
As previously, absorbing $Y_2$ in (\ref{5eq:264}$,c$) gives
\be
\label{5eq:266}
\|Y_2\|_{\sigma,\zeta} \,\leq\, C \frac{1}{\la\omega\ra^2 } \big[\|F\|_{\sigma,\zeta} +\|Y_1\|_{\sigma,\zeta} +|\zeta|C_2 I_\theta \big]
\ee
so that we obtain, thanks to (\ref{5eq:263}$,i$) in $(d)$:
\be
\label{5eq:267}
\ba{lrcl}
(a)&|\zeta p_0| &\leq& C\la\omega\ra^2 M_\sigma \big[\|F\|_{\sigma,\zeta} +\|Y_1\|_{\sigma,\zeta}  +|\zeta|\, C_2 \,I_\theta \big]  \medbreak\\
(b)&I_1 &\leq& C \la\omega\ra^2 M_\sigma M_{-\sigma-2} \big[\|F\|_{\sigma,\zeta} +\|Y_1\|_{\sigma,\zeta} +|\zeta|\, C_2\, I_\theta  \big]\medbreak\\
(c)&I_2 &\leq& C \big[\|F\|_{\sigma,\zeta} +\|Y_1\|_{\sigma,\zeta} +|\zeta|\, C_2\, I_\theta  \big] \medbreak\\
(d)& I_\theta &\leq& C \big[ \,\|F\|_{\sigma,\zeta} +  |\zeta|\, C_2\, \|Y_1\|_{\sigma,\zeta} +|\zeta|^2\, C_2^2\, I_\theta  + |\zeta|\, C_1\, I_1\,\big]
\ea
\ee
From corollary \ref{5cor:2}  we have $\la\omega\ra^2 M_\sigma M_{-\sigma}\leq \la\omega\ra^{-1}$  for $\sigma\in]-\frac32,-\frac12[$, thus we can absorb $Y_1$ in (\ref{5eq:267}$,b$) to obtain:
\be
\label{5eq:271}
\|Y_1\|_{\sigma,\zeta} \,\leq\, C   \la\omega\ra^2 M_\sigma M_{-\sigma}  \big[\|F\|_{\sigma,\zeta} +|\zeta|\,C_2\,I_\theta \big]
\ee
and
\be
\label{5eq:272}
\ba{lrcl}
(a)&|\zeta p_0| &\leq& C\la\omega\ra^2 M_\sigma \big[\|F\|_{\sigma,\zeta}  +|\zeta|\, C_2 \,I_\theta \big]  \medbreak\\
(b)&I_1 &\leq& C \la\omega\ra^2 M_\sigma M_{-\sigma-2} \big[\|F\|_{\sigma,\zeta}  +|\zeta|\, C_2\, I_\theta  \big]\medbreak\\
(c)&I_2 &\leq& C \big[\|F\|_{\sigma,\zeta}  +|\zeta|\, C_2\, I_\theta  \big]
\ea
\ee
Let us now examine $(d)$. Using lemma \ref{5lem:8} and corollary \ref{5cor:2}  we get:
\be
\label{5eq:273}
\ba{lrcl}
(d)& I_\theta &\leq& C \big[ \,\|F\|_{\sigma,\zeta} +  |\zeta|\, C_2\, \|Y_1\|_{\sigma,\zeta} +|\zeta|^2\, C_2^2\, I_\theta  + |\zeta|\, C_1\, I_1\,\big]  \medbreak\\
&&\leq& C \,\|F\|_{\sigma,\zeta} \,\big[ 1 +  |\zeta|\, C_2\, \la\omega\ra^2 M_\sigma M_{-\sigma} + |\zeta|\, C_1\,\la\omega\ra^2 M_\sigma M_{-\sigma-2} \big] + \smallbreak\\
&&  &  C\, |\zeta|\, C_2 \,I_\theta \,\big[ |\zeta|\, C_2 \,\la\omega\ra^2 \,M_\sigma M_{-\sigma}  + |\zeta|\, C_2 + |\zeta|\, C_1\ \, \la\omega\ra^2 \,M_\sigma M_{-\sigma-2} \,\big]  \medbreak\\
&&\leq& C  \,\big[ 1  + |\zeta|\, C_1\,\la\omega\ra^2 M_\sigma M_{-\sigma-2} \big] \,\big[ \, \|F\|_{\sigma,\zeta}+  \, |\zeta|\, C_2\,I_\theta \,\big] \medbreak\\
&&\leq& C  \,\big[ 1  + |\zeta|\, C_1\,\mbbm{1}_{\la\zeta\ra\geq\la\omega\ra^\kappa} + C_1\,\mbbm{1}_{\la\zeta\ra\leq\la\omega\ra^\kappa}  \big] \,\big[ \, \|F\|_{\sigma,\zeta}+  \, |\zeta|\, C_2\,I_\theta \,\big] \medbreak\\
&&\leq& C  \,\big[ \, \|F\|_{\sigma,\zeta}+  \, |\zeta|\, C_2\,I_\theta \,\big] 
\ea
\ee
According to corollary \ref{5cor:2}, $|\zeta|\,C_2$ tends to $0$ as $\la\omega\ra$ goes to infinity, we thus absorb $I_\theta$ in $(d)$ then in $(a,b,c)$ to conclude.\\
\eproof

\paragraph{Proof of lemma \ref{5lem:8}.} Let $s<\frac12$ and $\varphi \in H^{s+2}_\zeta$, with $\varphi = \sum_l \varphi_l e_l(z)$. Define $\phi= \int_0^z \varphi(z')\,dz'=\sum_k \phi_k e_k(z)$. Let us evaluate $\phi_k$:
\be
\label{5eq:247}
\ba{rcl}
\phi_k &=& C_0 \int_0^a \phi(z) \sin(\frac{k\pi z}a)\,dz \smallbreak \\
&=&C_0 \int_0^a \big[\int_0^z  \sum_l \varphi_l \sin(\frac{l\pi z'}a)\, dz' \big]\sin(\frac{k\pi z}a)\,dz \smallbreak \\
&=& C \int_0^a \big[\sum_l \varphi_l  \frac1l(1-\cos(\frac{l\pi z}a))\big] \sin(\frac{k\pi z}a)\,dz \smallbreak \\
&=&C \sum_l \int_0^a \frac{\varphi_l}l \big[ \sin(\frac{k\pi z}a) - \frac12 (\sin(\frac{(k+l)\pi z}a) + \sin(\frac{(k-l)\pi z}a)) \big] \,dz \smallbreak \\
&=& C \sum_{l}  \frac{\varphi_l}l \big[ \frac{1-(-1)^k}k - \frac{(1-(-1)^{k+l})k}{k^2-l^2} \big]
\ea
\ee
therefore
\be
\label{5eq:248}
|\phi_k |\,\sim \, C \sum_{l}  \frac{|\varphi_l|}l \Big[ \frac{1}k + \frac{k}{(k+l)(|k-l|+1)} \Big]
\ee
Let us first examine the case $s\in]-\frac12,\frac12[$. We have
\be
\label{5eq:257}
\ba{rcl}
\|\phi\|^2_s&\sim&\sum_k (1+k^2+\zeta^2)^s \bigg( \sum_{l}  \frac{|\varphi_l|}l \big[  \frac1k + \frac{k}{(k+l)(|k-l|+1)} \big] \bigg) ^2\\
&\leq& \|\varphi\|^2_s\sum_k (1+k^2+\zeta^2)^s \sum_{l} \frac1{l^2 (1+l^2+\zeta^2)^s } \big[ \frac1k + \frac{k}{(k+l)(|k-l|+1)} \big]^2
\ea
\ee
When $l$ is much smaller or much larger than $k$, $[ \frac1k + \frac{k}{(k+l)(|k-l|+1)} ]$ is of the same order as $\frac1k$; when $l$ and $k$ are of the same order, $[ \frac1k + \frac{k}{(k+l)(|k-l|+1)} ]$ is of the same order as $1$. Therefore we get:
\be
\label{5eq:250}
\ba{rcl}
\|\phi\|^2_s&\leq& C\,\|\varphi\|^2_s  \sum_k (1+k^2+\zeta^2)^s  \bigg( \frac1{k^2} \sum_{l} \big(\frac1{l^2 (1+l^2+\zeta^2)^s }\big) + \frac1{k^2 (1+k^2+\zeta^2)^s }\bigg)\smallbreak\\
&\leq& C\,\|\varphi\|^2_s \bigg( \sum_k \frac{(1+k^2+\zeta^2)^s }{k^2} \sum_{l} \frac1{l^2 (1+l^2+\zeta^2)^s } + \sum_k \frac1{k^2} \bigg)
\ea
\ee
For $s$ positive one has:
\be
\label{5eq:259}
\ba{rcl}
\sum_k \frac{(1+k^2+\zeta^2)^s }{k^2} &\sim& \int_1^\zeta \frac{(x^2+\zeta^2)^s }{x^2} dx + \int_\zeta^{+\infty} \frac{(x^2+\zeta^2)^s }{x^2} dx\smallbreak\\
&\sim &  \zeta^{2s} +  \zeta^{2s-1}\quad\sim\quad \zeta^{2s}
\ea
\ee
For $s$ positive, we have also:
\be
\label{5eq:260}
\ba{rcl}
\sum_{l} \frac1{l^2 (1+l^2+\zeta^2)^s } \quad\sim\quad \zeta^{-2s}
\ea
\ee
so finally one gets, for $s\in]-\frac12,\frac12[$:
\be
\label{5eq:260bis}
\ba{rcl}
\sum_{l} \frac{(1+l^2+\zeta^2)^s}{l^2  } \quad\sim\quad \zeta^{2s}
\ea
\ee
such that (\ref{5eq:246}) is established, for $s\in]-\frac12,\frac12[$.\\
We turn now to the case $s\in]-\frac32,-\frac12[$. Let us set
\be
\label{5eq:249}
\ba{rcl}
I_\varphi^2 &=& g_{\omega,\zeta}^2 \|\varphi\|^2_{s,\zeta} + \|\varphi\|^2_{s+2,\zeta}\smallbreak\\
&=& \sum_l (1+l^2+\zeta^2)^s(g_{\omega,\zeta}^2 + (1+l^2+\la\zeta\ra^2)^2) |\varphi_l|^2\smallbreak\\
&\sim& \sum_l (1+l^2+\zeta^2)^s((g_{\omega,\zeta}+\la\zeta\ra^2)^2 + l^4) |\varphi_l|^2\
\ea
\ee
where $g_{\omega,\zeta}$ can be either $\la\omega\ra^2$ or $\frac{M_{-s-2}}{M_{-s}}$.\\
Let us estimate $\|\phi\|_s$ as a function of  $I_\varphi$. From (\ref{5eq:248}) and the Cauchy-Schwarz inequality one has:
\be
\label{5eq:258}
\ba{l}
\|\phi\|^2_s\leq I_\varphi^2 \sum_k (1+k^2+\zeta^2)^s \sum_{l} \frac1{l^2 (1+l^2+\zeta^2)^s ((g_{\omega,\zeta}+\la\zeta\ra^2)^2 + l^4)} \big[ \frac1k + \frac{k}{(k+l)(|k-l|+1)} \big]^2
\ea
\ee
As previously and as for the estimation of $M_\sigma$, we have
\be
\label{5eq:252}
\ba{l}
\dsum_{l<< k\textrm{ or } l>>k} \frac1{l^2 (1+l^2+\zeta^2)^s ((g_{\omega,\zeta}+\la\zeta\ra^2)^2 + l^4)} \bigg[ \frac1k + \frac{k}{(k+l)(|k-l|+1)} \bigg]^2 \smallbreak\\
\qquad\qquad\qquad\leq\quad C \quad\dfrac1{k^2} \dsum_l \frac1{l^2 (1+l^2+\zeta^2)^s ((g_{\omega,\zeta}+\la\zeta\ra^2)^2 + l^4)}\smallbreak\\
\qquad\qquad\qquad\leq\quad C\quad \dfrac1{k^2} [ I_1+I_2+I_3 ]
\ea
\ee
with
\be
\label{5eq:253}
\ba{rcl}
I_1&=& \dint_1^\zeta \dfrac{dx}{x^2 (1+x^2+\zeta^2)^s ((g_{\omega,\zeta}+\la\zeta\ra^2)^2 + x^4)}\sim\dfrac{\zeta^{-2s}}{(g_{\omega,\zeta}+\la\zeta\ra^2)^2}\bigbreak\\
I_2&=& \dint_\zeta^{(g_{\omega,\zeta}+\la\zeta\ra^2)^{1/2}} \dfrac{dx}{x^2 (1+x^2+\zeta^2)^s ((g_{\omega,\zeta}+\la\zeta\ra^2)^2 + x^4)}\medbreak\\
&\sim &\dfrac{1}{(g_{\omega,\zeta}+\la\zeta\ra^2)^2} \dint_\zeta^{(g_{\omega,\zeta}+\la\zeta\ra^2)^{1/2}} \dfrac{dx}{x^{2s+2}}\quad\sim \quad \dfrac{(g_{\omega,\zeta}+\la\zeta\ra^2)^{-s-1/2}}{(g_{\omega,\zeta}+\la\zeta\ra^2)^2}\bigbreak\\
I_3&=& \dint_{(g_{\omega,\zeta}+\la\zeta\ra^2)^{1/2}}^{+\infty} \dfrac{dx}{x^2 (1+x^2+\zeta^2)^s ((g_{\omega,\zeta}+\la\zeta\ra^2)^2 + x^4)}\sim \dfrac1{(g_{\omega,\zeta}+\la\zeta\ra^2)^{s+5/2}}
\ea
\ee
Hence
\be
\label{5eq:254}
\ba{l}
\dsum_{l<< k\textrm{ or } l>>k} \frac1{l^2 (1+l^2+\zeta^2)^s (g_{\omega,\zeta}^2 + l^4+\zeta^4)} \bigg[ \frac1k + \frac{k}{(k+l)(|k-l|+1)} \bigg]^2 \smallbreak\\
\qquad\qquad\qquad\leq\qquad C \qquad\dfrac1{k^2} \dfrac{\zeta^{-2s}+  (g_{\omega,\zeta}+\la\zeta\ra^2)^{-s-1/2}  }{(g_{\omega,\zeta}+\la\zeta\ra^2)^2}
\ea
\ee
We have also:
\be
\label{5eq:255}
\sum_k (1+k^2+\zeta^2)^s \dfrac1{k^2} \sim \zeta^{2s}
\ee
Let us now examine (\ref{5eq:258}) for $l$ and $k$ of the same order:
\be
\label{5eq:256}
\ba{l}
\sum_k (1+k^2+\zeta^2)^s \dsum_{k,l\simeq k} \dfrac1{l^2 (1+l^2+\zeta^2)^s ((g_{\omega,\zeta}+\la\zeta\ra^2)^2 + l^4)} \bigg[ \frac1k + \frac{k}{(k+l)(|k-l|+1)} \bigg]^2\smallbreak\\
\quad\quad \sim \dsum_k (1+k^2+\zeta^2)^s \dfrac1{k^2 (1+k^2+\zeta^2)^s ((g_{\omega,\zeta}+\la\zeta\ra^2)^2+k^4)} \big[ \frac1k + C \big]^2\smallbreak\\
\quad\quad \sim \dsum_k  \dfrac1{k^2 ((g_{\omega,\zeta}+\la\zeta\ra^2)^2+k^4)} \smallbreak\\
\quad\quad \leq C \dfrac1{(g_{\omega,\zeta}+\la\zeta\ra^2)^2}
 \ea
\ee
We have proved that, for $s\in]-\frac32,-\frac12[$:
\be
\label{5eq:261}
\ba{rcl}
\|\phi\|_s \quad\leq\quad I_\varphi \quad \bigg(( g_{\omega,\zeta}+\la\zeta\ra^2)^{-1}    + \la\zeta\ra^{s}(g_{\omega,\zeta}+\la\zeta\ra^2)^{-s/2-5/4} \bigg)
\ea
\ee
Let us set
\be
\label{5eq:278}
\ba{rcl}
C_2(\omega,\zeta) &=& ( \la\omega\ra^2 +\la\zeta\ra^2)^{-1}   + \la\zeta\ra^{s}( \la\omega\ra^2+\la\zeta\ra^2)^{-s/2-5/4}\medbreak\\
C_1(\omega,\zeta) &=& ( \frac{M_{-s-2}}{M_{-s}} +\la\zeta\ra^2)^{-1}   + \la\zeta\ra^{s}( \frac{M_{-s-2}}{M_{-s}}+\la\zeta\ra^2)^{-s/2-5/4}
\ea
\ee
We have easily, considering separately the cases  $s\geq -1$ and $s\leq -1$:
\be
\label{5eq:270}
\ba{rcl}
|\zeta| \,C_2 &\sim &|\zeta|\, ( \la\omega\ra^2 +\la\zeta\ra^2)^{-1}   + \la\zeta\ra^{s+1}( \la\omega\ra^2+\la\zeta\ra^2)^{-s/2-5/4}\smallbreak\\
&\leq& C \la\omega\ra^{-1}    +  \la\zeta\ra^{s+1} \la\omega\ra^{-s-5/2} \smallbreak\\
&\leq& C \la\omega\ra^{-1} 
\ea
\ee
which gives in particular (\ref{5eq:263}$,i$).\\
Let us now turn to $C_1$, and set
\be
\label{5eq:279bis}
\kappa \,=\, \frac{2(-s-2)+1}{2(-s-2)}  \,=\, \frac{2s+3}{2s+4} 
\ee
For $s\in]-\frac32,-\frac12[$ we have:
\be
\label{5eq:279}
\ba{rcl}
 \frac{M_{-s-2}}{M_{-s}} &\sim& ( \la\zeta\ra^{s+2}\la\omega\ra^{-2} \mbbm{1}_{\la\zeta\ra\geq\la\omega\ra^\kappa} + \la\omega\ra^{s-1/2} \mbbm{1}_{\la\zeta\ra\leq\la\omega\ra^\kappa} ) \la\zeta\ra^{-s} \la\omega\ra^{2}\smallbreak\\
 &\sim&  \la\zeta\ra^{2} \mbbm{1}_{\la\zeta\ra\geq\la\omega\ra^\kappa} +  \la\zeta\ra^{-s}\la\omega\ra^{s+3/2} \mbbm{1}_{\la\zeta\ra\leq\la\omega\ra^\kappa}  
\ea
\ee
Thanks to (\ref{5eq:278}), we get:
\be
\label{5eq:280}
\ba{rcl}
\la\zeta\ra\, C_1 \mbbm{1}_{\la\zeta\ra\geq\la\omega\ra^\kappa} &\sim& (\la\zeta\ra^{-2}   + \la\zeta\ra^{-5/2})\,\la\zeta\ra \, \mbbm{1}_{\la\zeta\ra\geq\la\omega\ra^\kappa}
\ea
\ee
which gives (\ref{5eq:263}$,ii$), as $\kappa >0$. Then
\be
\label{5eq:281}
\ba{rcl}
 C_1 \mbbm{1}_{\la\zeta\ra\leq\la\omega\ra^\kappa} &\sim& ( \la\zeta\ra^{-s}\la\omega\ra^{s+3/2} +\la\zeta\ra^2)^{-1}   + \la\zeta\ra^{s}(\la\zeta\ra^{-s}\la\omega\ra^{s+3/2}  +\la\zeta\ra^2)^{-s/2-5/4}\smallbreak\\
&\leq& ( \la\omega\ra^{s+3/2} +\la\zeta\ra^2)^{-1}   + \la\zeta\ra^{s}(\la\zeta\ra^{-s}\la\omega\ra^{s+3/2}  +\la\zeta\ra^{-s})^{-s/2-5/4}\smallbreak\\
&\leq& ( \la\omega\ra^{s+3/2} +\la\zeta\ra^2)^{-1}   + \la\zeta\ra^{s^2/2+9s/4}(\la\omega\ra^{s+3/2}  +1)^{-s/2-5/4}
\ea
\ee
with $\la\zeta\ra^{-s}>1$ and $\la\zeta\ra^2>\la\zeta\ra^{-s}$. Finally we notice that
\be
\label{5eq:282}
\ba{l}
s+\frac32>0\quad ;\quad \frac12 s^2+\frac94s <-1\quad ;\quad -\frac 12 s -\frac54 <-\frac12
\ea
\ee
therefore $C_1 \mbbm{1}_{\la\zeta\ra\leq\la\omega\ra^\kappa}$ tends to 0 as $\la\omega\ra$ tends to infinity and the lemma is proven.\\
\eproof

\subsection{End of the proof}\label{thm1:5}
\label{5sec:fin}

We are now in a position to complete the proof, thanks to estimates (\ref{5eq:174}) of thm \ref{5prop:8}.\\
For $u$ and $v$ one gets:
\be
\label{5eq:213}
\ba{rcl}
 \|Y_1\|_{\sigma+2,\zeta} &\leq& C \la\omega\ra^2 M_\sigma M_{-\sigma-2} \|F\|_{\sigma,\zeta} \medbreak\\
 \|Y_2\|_{\sigma+2,\zeta} &\leq& C \|F\|_{\sigma,\zeta}
\ea
\ee
According to corollary \ref{5cor:2}, $\la\omega\ra^2 M_\sigma M_{-\sigma-2} \leq C$, therefore 
\be
\label{5eq:214}
 \|Y\|_{\sigma+2,\zeta} \leq C \|F\|_{\sigma,\zeta}
\ee
Inverting the Fourier-Laplace transform, one gets
\be
\label{5eq:215}
\ba{c}
Y \in L^2(\mbb{R}_+;H^{\sigma+2}_\zeta)\medbreak\\
 \|Y\|_{L^2(\mbb{R}_+;H^{\sigma+2}_\zeta)} \quad\leq\quad C \|F\|_{L^2(\mbb{R}_+;H^{\sigma}_\zeta)}
\ea
\ee
hence (\ref{5eq:29}) follows, for $u$ and $v$.\\
The same argument applies to  $\theta$:
\be
\label{5eq:216}
\ba{c}
\theta \in L^2(\mbb{R}_+;H^{\sigma+2}_\zeta)\medbreak\\
 \|\theta\|_{L^2(\mbb{R}_+;H^{\sigma+2}_\zeta)} \quad\leq\quad C \|F\|_{L^2(\mbb{R}_+;H^{\sigma}_\zeta)}
\ea
\ee
and (\ref{5eq:29}) is established for $\theta$, so that (\ref{eq:14}) and  (\ref{eq:15}) are proven.\\
For $\sigma \in ]-\frac12,\frac12[$, one has $M_\sigma \sim \frac{\la\zeta\ra^{-\sigma}}{\la\omega\ra^2}$ and with (\ref{5eq:174},a) one obtains
\be
\label{5eq:217}
\la\zeta\ra^{\sigma+1}|p_0| \quad\leq\quad C \|F\|_{\sigma,\zeta} 
\ee
then $p_0=q$ satisfies
\be
\label{5eq:218}
\ba{c}
q \in L^2(\mbb{R}_+;H^{\sigma+1}(\mbb{T}^2))\medbreak\\
\|q\|_{L^2(\mbb{R}_+;H^{\sigma+1}(\mbb{T}^2))} \quad\leq\quad C \|F\|_{L^2(\mbb{R}_+;\mc{H}^{\sigma})}
\ea
\ee
hence (\ref{5eq:219}) and (\ref{eq:18}) follow, for $\sigma \in ]-\frac12,\frac12[$.\\
For $\sigma \in ]-\frac32,-\frac12[$, we split $p_0=q$ in $q_1+q_2$ as follows:
\be
\label{5eq:220}
q_1=p_0 \,\mbbm{1}_{\la\zeta\ra\leq\la\omega\ra^\kappa},\quad q_2=p_0 \,\mbbm{1}_{\la\zeta\ra\geq\la\omega\ra^\kappa}
\ee
Similarly, lemma \ref{5prop:1} gives for  $q_2$:
\be
\label{5eq:221}
\la\zeta\ra^{\sigma+1}|q_2| \quad\leq\quad C \|F\|_{\sigma,\zeta} \\
\ee
thus
\be
\label{5eq:222}
\ba{c}
q_2 \in L^2(\mbb{R}_+;H^{\sigma+1}(\mbb{T}^2))\medbreak\\
\|q_2\|_{L^2(\mbb{R}_+;H^{\sigma+1}(\mbb{T}^2))} \quad\leq\quad C \|F\|_{L^2(\mbb{R}_+;\mc{H}^{\sigma})}
\ea
\ee
For $q_1$,  lemma \ref{5prop:1} gives:
\be
\label{5eq:223}
\la\omega\ra^{1/2+\sigma}\la\zeta\ra |q_1| \quad\leq\quad C \|F\|_{\sigma,\zeta} 
\ee
As
\be
\label{5eq:223bis}
\kappa = \frac{2\sigma+1}{2\sigma}\, \in\,\, ]0,\frac23[
\ee
we get 
\be
\label{5eq:223ter}
\la\zeta\ra \leq \la\omega\ra^\kappa =(|\tau|+\la\zeta\ra^2)^{1/2} \,\Rightarrow \, \la\zeta\ra \leq |\tau|^{1/2}
\ee
Therefore $\la\omega\ra \sim \la\tau\ra$ in the support of $q_1$, so (\ref{5eq:223}) becomes
\be
\label{5eq:225}
\ba{c}
q_1 \in H^{\sigma/2+1/4}(\mbb{R}_+;H^{1}(\mbb{T}^2))\medbreak\\
\|q_1\|_{H^{\sigma/2+1/4}(\mbb{R}_+;H^{1}(\mbb{T}^2))} \quad\leq\quad C \|F\|_{L^2(\mbb{R}_+;\mc{H}^{\sigma})}
\ea
\ee
\eproof

\section{Remarks and further results}

\subsection{An explicit formula for the pressure}\label{5rmq:12} 
Study of the uncoupled system gives the following formula for $p_0$:
\be
\label{5eq:226}
 p_0  =- \frac1a \frac{\omega^2}{\zeta^2} \Big[\mc{N}(\frac{\omega a}{\sqrt{\nu}})\Big]^{-1} \int_0^a  (\omega^2 -\nu \dzz)^{-1}[i\xi f_1+i\eta f_2]  \,dz
\ee
where $\mc{N}$ is defined by formula (\ref{5eq:147}). For the coupled system,  (\ref{5eq:198}) enables us to write
\be
\label{5eq:227}
\ba{l}
 p_0  =- \frac1a \frac{\omega^2}{\zeta^2} \big[\mc{N}(\frac{\omega a}{\sqrt{\nu}})\big]^{-1} \int_0^a  (\omega^2 -\nu \dzz)^{-1}[i\xi (f_1+\alpha v  -i\xi\beta\int_0^z \theta )\\
\quad\qquad\qquad\qquad\qquad\qquad\qquad\qquad\qquad +i\eta (f_2 -\alpha u -i\eta\beta \int_0^z \theta) ]  \,dz
\ea
\ee
According to theorem \ref{5prop:8}, the terms involving $u$ and $\theta$ are smooth enough that the singular term of the pressure is given by
\be
\label{5eq:228}
q_1(\tau,\xi,\eta)=-\,\mbbm{1}_{\la\zeta\ra\leq\la\omega\ra^\kappa}\, \frac1a \frac{\omega^2}{\zeta^2} \Big[\mc{N}(\frac{\omega a}{\sqrt{\nu}})\Big]^{-1} \int_0^a  (\omega^2 -\nu \dzz)^{-1}[i\xi f_1+i\eta f_2]  \,dz
\ee
with $\mc{N}$ defined by (\ref{5eq:147}).

\subsection{Proof of corollary \ref{cor:1}}\label{5sec:pvcor}

Under the assumptions of corollary \ref{cor:1},  $\varphi X$ satisfies equation (\ref{eq:13}) with right-hand side $\varphi F+\varphi'(t)X$ in $(L^2(0,T;\mc{H}^{-1}))^3$ and with support in $t> 0$. Theorem \ref{thm:1} gives (\ref{eq:23}). Thanks to the previous remark, we know that  $q$ splits in $q_1+q_2$, with $q_1$ less smooth than $q_2$, and we have explicitly (in Fourier variables)
\be
\label{5eq:275}
q_1(\tau,\xi,\eta)=-\,\mbbm{1}_{\la\zeta\ra\leq\la\omega\ra^\kappa}\, \frac1a \frac{\omega^2}{\zeta^2} \Big[\mc{N}(\frac{\omega a}{\sqrt{\nu}})\Big]^{-1} \int_0^a  (\omega^2 -\nu \dzz)^{-1}[i\xi (\varphi f_1)_{\tau,\zeta}+i\eta (\varphi f_2)_{\tau,\zeta}]  \,dz
\ee
(indeed, the contribution of $\varphi'(t)X$ is smooth, because $C\in L^2(0,T;\mc{V})$)\\
Thus the following equivalence holds:
\be
\label{5eq:276}
\ba{rcl}
&&q_1 \in L^2(0,T;L^2(\mbb{T}^2))\smallbreak\\
 &\Leftrightarrow& \frac1{\zeta^2} \int_0^a  (\omega^2 -\nu \dzz)^{-1}[i\xi (\varphi f_1)_{\tau,\zeta}+i\eta (\varphi f_2)_{\tau,\zeta}]  \,dz \in L^2(\tau;\ell^2_\zeta) \smallbreak\\
 &\Leftrightarrow& \Delta_2^{-1} \big[ \int_0^a (\dt-\nu\Delta)^{-1}[\varphi\dx f_1+\varphi\dy f_2]\,dz \big] \in L^2(0,T;L^2(\mbb{T}^2))
\ea
\ee
\eproof

\subsection{A counter-example to maximal estimates}\label{5rmq:13} 
In this paragraph, we construct a counter-example to the maximal estimate  (\ref{eq:21}) for $\sigma<-1/2$. 
It is sufficient to find $F\in (L^2(\mbb{R},\mc{H}^\sigma))^3$ such that the gradient of the associated $p_0$ is not in $L^2$ in time, as each other term of the equations (\ref{eq:13})(except $\dt X$) are in $L^2$ in time. \\
Let $\sigma\in]-\frac32,-\frac12[$ and let $\alpha$ be such that $\alpha \in ]\sigma+\frac12,0[$. Let $g(t)\in L^2(\mbb{R})$ with support in $t> 0$.
Let $f\in \mc{H}^\sigma$, independent of time, be defined by
\be
\label{5eq:229}
f(x,y,z) = \sum_k k^{-\alpha}\,e_k(z) e^{ix}
\ee
If $\zeta=(1,0)$ then $f_{\zeta,k}=k^{-\alpha}$ and if $\zeta\neq(1,0)$ then  $f_{\zeta,k}=0$.\\
Now let  $F \in (L^2(\mbb{R},\mc{H}^\sigma))^3$ with support in $t>0$ defined by  $F=(fg,0,0)$ and let $p$ be the pressure, solution of (\ref{eq:13}) with right hand side $F$. We write $p=q_1+p_1$ where $\nabla p_1$ is as smooth as $F$ is, and $q_1$ is explicitly given by formula  (\ref{5eq:228}) (omitting the high frequency cut-off):
\be
\label{5eq:230}
q_1(\tau,\xi,\eta)=- \frac1a \frac{\omega^2}{\zeta^2} \Big[\mc{N}(\frac{\omega a}{\sqrt{\nu}})\Big]^{-1} \int_0^a  (\omega^2 -\nu \dzz)^{-1}[i\xi f(\xi,\eta,z) g(\tau)]  \,dz
\ee
Therefore $q_1(t,x,y)=q_1(t)\, e^{ix}$ and the Fourier transform of $q_1(t)$ is
\be
\label{5eq:231}
\ba{lrcl}
&q_1(\tau)&=& g(\tau)\, m(\tau)\medbreak\\
with&m(\tau) &=&  \frac1a (\tau -i ) \Big[\mc{N}(\frac{\omega a}{\sqrt{\nu}})\Big]^{-1} \int_0^a  (\omega^2 -\nu \dzz)^{-1}[ f_{(1,0)}(z)  ]  \,dz
\ea
\ee
with $\omega^2=i\tau+\zeta^2$. Let us now calculate $m(\tau)$:
\be
\label{5eq:286}
\ba{rcl}
m(\tau) &=&  \frac1a (\tau - i) \Big[\mc{N}(\frac{\omega a}{\sqrt{\nu}})\Big]^{-1} \dint_0^a  \dsum_{k} \dfrac{f_{k}}{\omega^2+\frac{\nu k^2 \pi^2}{a^2}} e_k(z) \,dz\smallbreak\\
&=&  \frac{C_a}a (\tau - i) \Big[\mc{N}(\frac{\omega a}{\sqrt{\nu}})\Big]^{-1} \dsum_{k\textrm{ odd}}\dfrac{k^{-\alpha}}{k\,(1+i\tau+\frac{\nu k^2 \pi^2}{a^2})}
\ea
\ee
In order to get a lower bound on  $|m(\tau)|$, let us define
\be
\label{5eq:287}
\ba{rcl}
S_k &=& \dsum_{k\textrm{ odd}}\dfrac{k^{-\alpha}}{k\,(1+i\tau+\frac{\nu k^2 \pi^2}{a^2})} \smallbreak\\
&=& \dsum_{k\textrm{ odd}}\dfrac{k^{-\alpha}(1+\frac{\nu k^2 \pi^2}{a^2})}{k\,((1+\frac{\nu k^2 \pi^2}{a^2})^2+\tau^2)} -i \dsum_{k\textrm{ odd}}\dfrac{k^{-\alpha} \tau }{k\,((1+\frac{\nu k^2 \pi^2}{a^2})^2+\tau^2)}
\ea
\ee
We have, as $\alpha+3>1$:
\be
\label{5eq:232}
\ba{rcl}
|(\tau-i)S_k|&\geq& |\Re ((\tau-i)S_k)| \,=\,| \tau \Re(S_k) + \Im(S_k)\smallbreak|\\
&=&|\tau| \dsum_{k\textrm{ odd}}\dfrac{ k^{-\alpha}}{k\,((1+\frac{\nu k^2 \pi^2}{a^2})^2+\tau^2)} (1+\frac{\nu k^2 \pi^2}{a^2}-1)\smallbreak\\
&=&C|\tau| \dsum_{k\textrm{ odd}}\dfrac{ k^{2-\alpha} }{k\,((1+\frac{\nu k^2 \pi^2}{a^2})^2+\tau^2)} \smallbreak\\
&\geq&C|\tau| \dsum_{k \geq |\tau|^{1/2}}\dfrac{ k^{2-\alpha} }{k\,((1+\frac{\nu k^2 \pi^2}{a^2})^2+\tau^2)} \smallbreak\\
&\sim& |\tau| \dint_{|\tau|^{1/2}} \dfrac{x^{2-\alpha}\, dx}{x(x^4+\tau^2)}\smallbreak\\
&\sim&  |\tau|^{-\alpha/2}
\ea
\ee
Moreover we know that $\mc{N}^{-1}\rightarrow 1$ as $|\tau| \rightarrow +\infty$, thus we have, if $|\tau|$ is large enough
\be
\label{5eq:233}
|m(\tau)| \geq C  |\tau|^{-\alpha/2}
\ee
So finally, for  $|\tau|$ large enough
\be
\label{5eq:288}
|q_1(\tau)|\,\geq \,C\, |g(\tau)|\, |\tau|^{-\alpha/2}
\ee
Choose now $g\in L^2(\mbb{R})$, with support in $t>0$, such that $g\notin H^s(\mbb{R})$ for all $s>0$. Then $|g(\tau)|\, |\tau|^{-\alpha/2}$ is not in $L^2(\tau\in\mbb{R})$ (because $-\alpha>0$), and the pressure $q_1(t)e^{ix}$ is not $L^2$ in time and neither is its gradient.

\nocite{Nodet05}

\bibliographystyle{plain}
\bibliography{Biblio}

%
%
%
%
%
%
%
%
%
%
%
%
%
%
%
%
%
%
%

\end{document}